\documentclass{amsart}
\usepackage{amssymb}
\usepackage[latin1]{inputenc}
\usepackage[frenchb,germanb,dutch,italian,british]{babel}[1999/09/09]
\usepackage[all]{xy}

\makeatletter

\AtBeginDocument{% counteract harmful effects of Babel package
  \catcode`\!=12
  \catcode`\:=12
  \catcode`\;=12
  \catcode`\?=12
}

\renewcommand{\subsection}{\@startsection
	{subsection}%
	{2}%
	{0mm}%
	{1\baselineskip plus 0.25\baselineskip minus 0.25\baselineskip}%
	{-1em}%
	{\normalfont\normalsize\scshape}%
}

\newtheorem{thm}{Theorem}[section]
\newtheorem{prop}[thm]{Proposition}
\newtheorem{lem}[thm]{Lemma}
\newtheorem{cor}[thm]{Corollary}
\theoremstyle{definition}
\newtheorem{defn}[thm]{Definition}
\newtheorem{example}[thm]{Example}
\newtheorem{rem}[thm]{Remark}
\numberwithin{figure}{section}

\newcommand*{\bref}[1]{(\ref{#1})}
\newcommand*{\define}[2][]{\def\sequence@tmp{#1}%
  \ifx\sequence@tmp\@empty\def\sequence@tmp{#2}\fi
  \expandafter\index\expandafter{\sequence@tmp}\emph{#2}}
\providecommand{\url}[1]{\texttt{#1}}

\newcommand*{\closure}[2][]{\overline{#2}^{\scriptscriptstyle#1\relax}}
\newcommand*{\reciprocal}[1]{\frac{1}{#1}}
\newcommand*{\st}{\;\mathchoice{\vrule width .07em\relax}{\vert}{\vert}{\vert}\;}
\newcommand*{\contains}{\supseteq} 
\newcommand*{\sub}{\subseteq}
\newcommand*{\lee}{\leqslant} 
\newcommand*{\gee}{\geqslant} 
\newcommand*{\setprod}{\times}
\newcommand*{\defeq}{\mathrel{:=}} 
\newcommand*{\maps}{\colon}
\newcommand*{\union}{\cup} 
\newcommand*{\intersect}{\cap}
\newcommand*{\Intersection}{\bigcap} 
\newcommand*{\Union}{\bigcup}
\newcommand*{\Dirsum}{\bigoplus}
 
\newcommand*{\real}{\mathbb{R}}
\newcommand*{\Z}{\mathbb{Z}} 

\newcommand*{\nats}{\mathbb{N}} 
\newcommand*{\I}{\mathbb{I}}
\newcommand*{\continuum}{\mathfrak{c}} 
\newcommand*{\cardm}{\mathfrak{m}}
\newcommand*{\Bee}{\mathcal{B}} 
\newcommand*{\Eee}{\mathcal{E}}
\newcommand*{\Ess}{\mathcal{S}} 

\newcommand*{\Dee}{\mathcal{D}}

\renewcommand*{\emptyset}{\varnothing}

\newcommand*{\sequence@parenconstruct}[3]{\bgroup
	\mathchoice
	{\left#1#3\right#2}{#1#3#2}{#1#3#2}{#1#3#2}
\egroup}
\DeclareRobustCommand{\abs}{\sequence@parenconstruct{\lvert}{\rvert}}
\DeclareRobustCommand{\set}{\sequence@parenconstruct{\lbrace}{\rbrace}}
\def\seq(#1){\sequence@parenconstruct{(}{)}{#1}}

\newcommand*{\citet}[2][]{\def\sequence@tmpa{#1}% poor substitute for natbib
	\ifx\sequence@tmpa\@empty\cite{#2}\else
	\cite{#2}\space(#1)\fi}
\newcommand*{\citep}{\citet}

\makeatother

\begin{document}
\title{Sequential Convergence in Topological Spaces}
\author{Anthony Goreham}
\thanks{This survey was written as a dissertation for the Final Honour 
School of Mathematics, The Queen's College, Oxford University.}
\date{April 2001}
\begin{abstract}
In this survey, my aim has been to discuss the use of sequences and
countable sets in general topology. In this way I have been led to
consider five different classes of topological spaces: first countable
spaces, sequential spaces, Fréchet spaces, spaces of countable tightness
and perfect spaces. We are going to look at how these classes are
related, and how well the various properties behave under certain
operations, such as taking subspaces, products, and images under proper
mappings. Where they are \emph{not} well behaved we take the opportunity
to consider some relevant examples, which are often of special interest. 
For instance, we examine an example of a Fréchet space with unique
sequential limits that is not Hausdorff. I asked the question of whether
there exists in ZFC an example of a perfectly normal space that does not
have countable tightness: such an example was supplied and appears below. 
In our discussion we shall report two independence theorems, one of which
forms the solution to the Moore-Mrówka problem. The results that we prove
below include characterisation theorems of sequential spaces and Fréchet
spaces in terms of appropriate classes of continuous mappings, and the
theorem that every perfectly regular countably compact space has countable
tightness.
\end{abstract}
\maketitle

%\tableofcontents

%\clearpage 
%\pagenumbering{arabic}

%%%%%%%%%%%%%%%%%%%%%%%%%%%%%%%%%%%%%%%%%%%%%%%%%%%%%%%%%%%%%%%%%%%%%
%%%%%%%%%%%%%%%%%%%%%%%%%%%%%%%%%%%%%%%%%%%%%%%%%%%%%%%%%%%%%%%%%%%%%
%%%%                                                             %%%%
%%%%                      CHAPTER ONE                            %%%%
%%%%                                                             %%%%
%%%%%%%%%%%%%%%%%%%%%%%%%%%%%%%%%%%%%%%%%%%%%%%%%%%%%%%%%%%%%%%%%%%%%
%%%%%%%%%%%%%%%%%%%%%%%%%%%%%%%%%%%%%%%%%%%%%%%%%%%%%%%%%%%%%%%%%%%%%

\section{Introduction}\label{chap:intro}

\subsection*{Motivation}

For general topology, sequences are inadequate: everyone says so
\citep{Cartan37, Tukey40, Birkhoff37, Kelley55, Willard70, Engelking89}.
In Example~\ref{eg:betanats} we shall see some evidence to support this
claim. Undaunted, in the rest of this section I am going to set out how we
shall proceed in our discussion of the rôle of sequences, and countable
sets, in general topology. In Section~\ref{chap:firstcount} we review the
properties of sequences in the class of first countable spaces, most of
which are familiar. The work begins in earnest in
Section~\ref{chap:sequential} and Section~\ref{chap:frechet}, when we
consider sequential spaces and Fréchet spaces. It is easily shown that a
first countable space is Hausdorff if and only if its sequential limits
are unique. However, this is no longer true for Fréchet spaces (even
compact Fréchet spaces), as we shall see in
Example~\ref{eg:frechet-hausdorff}. We are also going to prove
characterisation theorems of the classes of sequential spaces and Fréchet
spaces, in terms of appropriate continuous mappings (see
Theorems~\ref{thm:sequential-characterisation},~\ref{thm:frechet-characterisation}).

It is natural to ask what changes if, instead of studying the
sequences in a space, one studies its countable subsets. This
question leads us to consider the class of spaces of countable
tightness. The relationship between the properties of having
countable tightness and being a sequential space is of great
interest, being the issue at the centre of the famous Moore-Mrówka
problem. In one direction, it is a triviality to prove that every
sequential space has countable tightness. The question of whether
the converse is true for compact Hausdorff spaces occupied
mathematicians for over twenty years (see
Section~\ref{chap:countable-tightness}).

The class of perfect spaces is not often studied in this context, and it
does not immediately fit into the hierarchy formed by the others that we
have named (see Figure~\ref{fig:relations}). That a first countable normal
space need not be perfectly normal is well known. Furthermore,
\citet{Arens50} has given an example of a perfectly normal space that is
not sequential. It seems, however, that the question of whether a perfect
space need have countable tightness has never been properly considered.
This may be partly because the cardinal invariant~$\Psi$ associated with
the perfect property is rarely used (it does not even have a name; see
page~\pageref{perfectionindex}). A pertinent example was found by
Mr.~J.~Lo, and it appears as Example~\ref{eg:perfect-counttight}.

When suitable compactness conditions are imposed the results are more
encouraging. Indeed, we are going to prove that every perfectly regular
countably compact space has countable tightness
(Theorem~\ref{thm:perfectcompact-counttight}). Like the other classes that
we have defined, the class of perfect spaces contains all metric spaces.
So I asked what follows if a space, like a metric space, is both perfectly
normal and first countable. Of particular interest is the question of
whether compactness and sequential compactness are equivalent for such
spaces. The fact that this cannot be decided in ZFC is easily deduced from
known results (see Corollary~\ref{cor:perfectfirstcount-compact}).

\subsection*{An Example}

It is a familiar and useful fact that, in metric spaces, the crucial 
notions of closure, compactness and continuity can be described in terms 
of convergence of sequences. As a sample of the results available in 
metric spaces, we recall the following:
\begin{enumerate}
\item\label{intro-item:1} 
{\itshape The closure of a set consists of the limits of all sequences in 
that set;}
\item\label{intro-item:2} 
{\itshape The space is compact if and only if it is sequentially compact;}
\item\label{intro-item:3} 
{\itshape A function from the space into a topological space is continuous 
if and only if it preserves limits of sequences.}
\end{enumerate}
In general topological spaces, these results are no longer true,
as the following example shows.

\begin{example}\label{eg:betanats}
\index{compact space}\index{sequentially compact space} 
Let $\beta\nats$~denote the Stone-{\v C}ech compactification of the
natural numbers. Then the only sequences in that converge in~$\beta\nats$
are those which are eventually constant, i.e., they always take the same
value from some point on---see \citet[Cor.~3.6.15]{Engelking89}. We see
immediately that $\beta\nats$~is a compact space that is not sequentially
compact, which is a counter-example to item~\bref{intro-item:2} above. Now
let $x$~be any point in $\beta\nats\setminus\nats$. Since $\nats$~is dense
in~$\beta\nats$, we see that $x$ belongs to $\closure\nats$, the closure
of~$\nats$. However, no sequence in~$\nats$ converges to~$x$. This shows
that statement~\bref{intro-item:1} does not hold in~$\beta\nats$. Finally,
let $Y$~be the set of all points of~$\beta\nats$ in the discrete topology.
Then the topology on~$Y$ is strictly finer than that on~$\beta\nats$, so
the identity mapping $\iota\maps\beta\nats\to Y$ is \emph{not} continuous.
However, since the only convergent sequences in~$\beta\nats$ are the
trivial ones, we see that $\iota$~does preserve limits of sequences. This
shows that~\bref{intro-item:3} is not generally valid.
\qed
\end{example}

One can give simpler examples of spaces in
which~\bref{intro-item:1},~\bref{intro-item:3} do not hold; see
for example \citet[Sec.~3.1]{McCluskey97}. However,
$\beta\nats$~is the canonical counter-example
to~\bref{intro-item:2}. We shall return to this example of a space
that is not discrete, but in which there are no non-trivial
convergent sequences.

There will be more examples in subsequent sections of spaces in
which sequences fail to describe important topological properties
of the space. Thus it would appear that sequences are inadequate
to the task of describing topology, in general. Why should we be
worried by this? After all, there are many properties of metric
spaces that fail to generalise to topological spaces; this is what
makes topological spaces so interesting. The difficulty is that
convergence of sequences is one of the most basic, and one of the
most intuitive, aspects of topology. The first attempts by
\citet{Frechet06} to find a general class of spaces suitable for
study of topological notions were centred around the convergence
of sequences. The fact that open sets, closure, continuity and
compactness cannot in general be described by sequential
convergence is thus a great inconvenience.

How do we respond to this? There is more than one possibility. One
approach is to stop thinking about convergence in terms of
sequences, and to substitute some more general theory of
convergence that is sufficient in any topological space. This is
done in the twin theories of nets and of filters, in which a
generalised convergence is developed that does describe open sets,
closure, continuity and compactness in a satisfactory manner; see
\citet[Chap.~4]{Willard70}. Another possibility is to accept that
sequences cannot always do everything, and to look to see what
they \emph{can} do, and when. This is the route that we shall
follow.

\subsection*{Manifesto}

We use the axiom of choice without question and without comment.
The set of all natural numbers is denoted~$\nats$; we do not
regard~$0$ as a natural number. A
\index{compact space}%
compact space is not necessarily Hausdorff, and in this regard our
notion of compactness differs from that of \citet{Bourbaki66}.
Likewise, we do not assume that a
\index{paracompact space}%
paracompact space is Hausdorff, nor that a Lindelöf space is
regular. A regular or a normal space is not necessarily~$T_1$;
thus $T_3$~is a condition strictly stronger than regularity, and
$T_4$~is strictly stronger than normality, for in both of these we
insist upon~$T_1$. A neighbourhood~$N$ of a point~$x$ need not
necessarily be an open set; all that is required is that $x$~be an
interior point of~$N$. We use the term `subsequence' in the
traditional sense, rather than in the sense of \citet{Kelley55}.
Other terms occurring are either defined as they are used, or are
thought to be unambiguous; in the latter case, the reader is
referred to \citet{Engelking89} for a definition.

%\medskip

\begin{figure}[h]
\begin{center}
\small
%% DIAGRAM: needs xymatrix package installed \usepackage[all]{xy}
\begin{displaymath}
\renewcommand{\objectstyle}{\text}
\renewcommand{\labelstyle}{\textstyle}
\xymatrix{%
metric space			\ar@{=>}[d]	\ar[r]				&\makebox[2in][l]{compactness $\iff$ sequential compactness}	\\
first countable space		\ar@{=>}[d]	\ar[r]	\ar[ur]|-{\backslash}	&\makebox[2in][l]{Hausdorff $\iff$ unique limits}		\\
Fréchet space			\ar@{=>}[d]	\ar[r]	\ar[ur]|-{\backslash}	&\makebox[2in][l]{closure $=$ set of limits}			\\
sequential space		\ar@{=>}[d]	\ar[r]	\ar[ur]|-{\backslash}	&\makebox[2in][l]{continuity $\iff$ preservation of limits}	\\
space of countable tightness				\ar[ur]|-{\backslash}
}
\end{displaymath}
\caption{The relations between some classes of spaces, and their properties.}
    \label{fig:relations}%
    \index{countable tightness}%
    \index{continuous mapping}%
    \index{Hausdorff property}%
    \index{sequential space}%
    \index{Fréchet space}%
    \index{first axiom of countability}%
    \index{compact space}%
    \index{sequentially compact space}%
\end{center}
\end{figure}
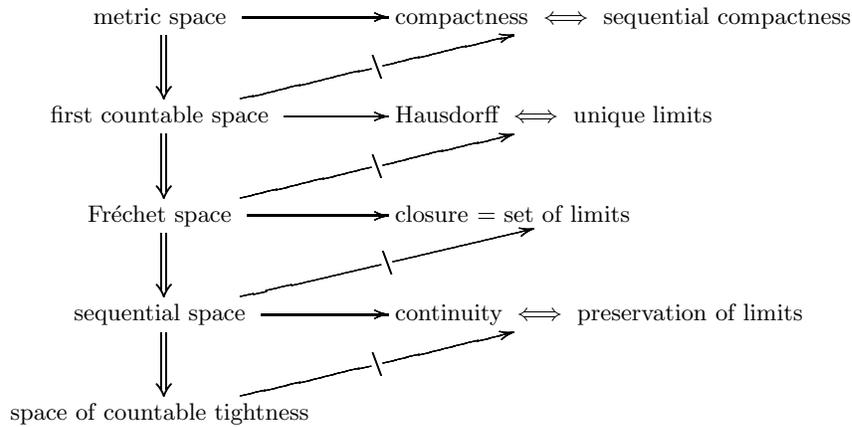

We shall examine, principally, four classes of spaces. Each class
is defined either directly in terms of sequences or in terms of
some sort of countability. These classes form a hierarchy in which
the class of spaces under consideration becomes wider, while the
results available become weaker. For an overview of this
hierarchy, see Figure~\ref{fig:relations}.

This diagram also indicates some of the properties of sequential
convergence in each class of spaces---the results and
counter-examples required to prove the various assertions in this
diagram will appear in later sections.

In each section we are looking at the sequential convergence
properties of the class of spaces under consideration. We also
consider the standard questions: `Is this class closed under
taking products, or subspaces, or (inverse) images under
appropriate continuous mappings?' In fact, the most useful
continuous mappings in this regard are proper mappings, insofar as
they are included in most of the results that we find. The answer
to these questions is frequently `no'---see
Figure~\ref{fig:class-properties}. This unfortunate fact nonetheless
gives us the opportunity to consider some interesting examples.

\begin{figure}[h]
\begin{center}
\small
\begin{tabular}{lcccc}
			&		&		&\multicolumn{2}{c}{Proper mappings:}	\\
Class of spaces		&Hereditary?	&Productive?	&Invariant?	&Inverse invariant?	\\
\hline\\[-1.5ex]
First countable		&Yes		&No		&No		&No			\\
Fréchet			&Yes		&No		&Yes		&No			\\
Sequential		&No		&No		&Yes		&No			\\
Countable tightness	&Yes		&No		&Yes		&No			\\
Perfect			&Yes		&No		&Yes		&No			\\
\end{tabular}
    \caption{}
    \label{fig:class-properties}%
    \index{proper mapping}%
    \index{first axiom of countability}%
    \index{Fréchet space}%
    \index{sequential space}%
    \index{countable tightness}%
    \index{perfect space}%
\end{center}
\end{figure}

Concerning the characterisation of classes of spaces in terms of
continuous images of metric spaces, the reader is referred to
Figure~\ref{fig:characterisations}. These results were achieved in
the 1960s by \citet{Ponomarev60}, \citet{Arhangelskij63} and
\citet{Franklin65}.

%\begin{figure}[tp]
\begin{figure}[h]
\begin{center}
\small
%%% DIAGRAM: needs xymatrix package installed \usepackage[all]{xy}
%\renewcommand{\objectstyle}{\text}
%    \begin{tabular}{ccc}
%      Space &
%      $+$ Condition &
%      $=$ Image of a metric space under \\
%      \hline\\[-1ex]
%      \xymatrix{
%        first countable space \ar @{=>} [d] \\
%        Fréchet space \ar @{=>} [d] \\
%        sequential space \\
%        } &
%      \xymatrix{
%        $T_0$ \\
%        $T_2$ \\
%        \\
%        } &
%      \xymatrix{
%        open onto mapping \ar @{=>} [d] \\
%        pseudo-open mapping \ar @{=>} [d] \\
%        quotient mapping \\
%        }
%    \end{tabular}

\begin{tabular}{ccc}
Space			&$+$ Condition	&$=$ Image of a metric space under	\\
\hline\\[-1.5ex]
first countable space	&$T_0$		&open onto mapping			\\
$\Downarrow$		&		&$\Downarrow$				\\
Fréchet space		&$T_2$		&pseudo-open mapping			\\
$\Downarrow$		&		&$\Downarrow$				\\
sequential space	&		&quotient mapping			\\
\end{tabular}

    \caption{Characterisations of spaces in terms of continuous mappings.}
    \label{fig:characterisations}%
    \index{pseudo-open mapping}%
    \index{open mapping}%
    \index{quotient mapping}%
    \index{sequential space}%
    \index{Fréchet space}%
    \index{first axiom of countability}%
  \end{center}
\end{figure}

Were I to include the proofs to all the relevant results mentioned above,
then this survey would be twice as long. Therefore I have carefully
selected which results and examples to consider in detail. I have been
particularly concerned for variety, so where two or more examples or
results use similar ideas, only one proof will be given. Moreover, very
few of the lemmas are proved: in most cases the reader can supply a proof
using the first idea that comes to mind.

%\pagebreak[3]
\subsection*{Definitions}

How do we interpret the statement that `sequences describe the
topology' in a space~$X$? Our usual definition of `topology' is
that the topology on~$X$ is the family of open subsets of~$X$.
Hence the most natural interpretation of our statement is that we
can give a necessary and sufficient condition for openness of a
subset of~$X$ in terms of sequences. The following is a simple
property of metric spaces:
\begin{enumerate}
  \setcounter{enumi}{3}
\item\label{intro-item:4} {\itshape A subset~$A$ of\/~$X$ is
    open if and only if whenever $\seq(x_i)$ is a sequence that
    converges to a point of~$A$, then $A$~contains all but finitely
    many of the members $x_i$ of the sequence.}
\end{enumerate}
We shall define the \emph{sequential} spaces to be the spaces in
which this statement holds. However, there are many other ways to
define a topology on~$X$, other than specifying the open sets. We
could equally well specify the closed sets in~$X$, or we could
define a Kuratowski closure operation on~$X$. In terms of closed
sets, we note that metric spaces have the following property:
\begin{enumerate}
  \setcounter{enumi}{4}
\item\label{intro-item:5} {\itshape A subset~$A$ of\/~$X$ is
    closed if and only if whenever $\seq(x_i)$ is a sequence in~$A$
    that converges to a point $x\in X$, then $x$~belongs to~$A$.}
\end{enumerate}
We shall see that properties~\bref{intro-item:4}
and~\bref{intro-item:5} are equivalent. It is perhaps surprising,
therefore, that they are \emph{not} equivalent to
property~\bref{intro-item:1}, which relates sequences to the
closure operator---see page~\pageref{intro-item:1}. It turns out
(see~Example~\ref{eg:sequential-frechet}) that
property~\bref{intro-item:1} is \emph{strictly} stronger than
properties \bref{intro-item:4}~and~\bref{intro-item:5}. Those
spaces that satisfy~\bref{intro-item:1} are the \emph{Fréchet}
spaces.

%\medskip

At this point we are going to define each of the main classes of spaces
that will be considered in future sections. I have chosen to do this,
rather than to withhold the definition of each class until we are ready to
study it, in order to facilitate appropriate comparison of these classes
as we go along.

To start with, it is well known that spaces satisfying the so-called first
axiom of countability share most of the nice properties, with regard to
sequences, of metric spaces.

\begin{defn}[first axiom of countability]\label{defn:firstcount}
Let $X$~be a topological space. We say that $X$~is \define[first countable 
space|see{first axiom of countability}]{first countable}, or that 
$X$~satisfies the \define{first axiom of countability}, if for every 
point~$x$ of\/~$X$ the neighbourhood system at~$x$ has a countable base.
\end{defn}

It is clear that every metric space is first countable: a countable 
neighbourhood base at the point~$x$ is the 
set~$\set{S(x,1/i)\st i\in\nats}$ of open spheres of radius~$1/i$ and 
centre~$x$, as $i$~runs through the natural numbers.

\begin{defn}[Fréchet space]\label{defn:frechet}
Let $X$~be a topological space. We say that $X$~is a 
\define[Fréchet space]{Fréchet} space, or that it has the
\define[Fréchet property|see{Fréchet space}]{Fréchet property}, 
if for every subset~$A$ of~$X$, a point $x\in X$ belongs to $\closure A$ 
if and only if there exists a sequence in~$A$ that converges to~$x$. 
Fréchet spaces are also known as 
\define[Fréchet-Urysohn space|see{Fréchet space}]{Fréchet-Urysohn} spaces.
\end{defn}

In analysis one often uses the fact that metric spaces are Fréchet. 
More generally, any first countable space is Fréchet:

\begin{lem}\label{lem:firstcount-frechet}
Let $X$~be a first countable space. Then $X$~is a Fréchet space.\qed
\end{lem}

The proof is easy to find; it appears in \citet[Thm.~3.10]{McCluskey97}, 
for instance.

It is convenient to adopt some terminology of \citet{Franklin65}.

\begin{defn}[sequentially open]
  Let $X$~be a topological space, and let $A$~be a subset of~$X$. We
  say that~$A$ is \define{sequentially open} if whenever $\seq(x_i)$
  is a sequence in~$X$ that converges to a point in~$A$, then there
  exists a natural number~$i_0$ such that $x_i$ belongs to~$A$
  whenever $i\gee i_0$. (In these circumstances we often say that
  $A$~contains \define{almost all} of the $x_i$, or that
  $\seq(x_i)$~is \define{eventually} in~$A$.) Similarly, we say that
  $A$~is \define{sequentially closed} if it contains the limits of all
  its sequences, i.e., if whenever there is a sequence in~$A$ that
  converges to some point~$x\in X$, then $x$~belongs to~$A$.
\end{defn}

From the definition of sequential convergence in a topological space,
every open set is sequentially open. Similarly, it is a trivial matter to
prove that every closed set is sequentially closed. That the converses are
not true in general spaces we can see by again considering $\beta\nats$,
which is non-discrete but might be said to be `sequentially discrete'.
Note that a subset~$A$ of~$X$ is sequentially open iff $X\setminus A$ is
sequentially closed. With this in mind, the next result is obvious.

\begin{lem}\label{lem:sequential-equiv}
Let $X$~be a topological space. Then the following are equivalent:
\begin{enumerate}
\item every sequentially open subset of~$X$ is open;
\item every sequentially closed subset of~$X$ is closed.\qed
\end{enumerate}
\end{lem}

\begin{defn}[sequential space]\label{defn:sequential}
Let $X$~be a topological space. We say that $X$~is a
\define[sequential space]{sequential} space if every sequentially open 
subset of~$X$ is open. By Lemma~\ref{lem:sequential-equiv} this is the 
same as saying that every sequentially closed set is closed.
\end{defn}

\begin{rem}[{\citet[p.~44]{Arhangelskij84a}}]\label{rem:frechetsequential}
What is the difference between the present property and that of being a
Fréchet space? We can see this more clearly by rephrasing the definition
in a third way: a space~$X$ is sequential if for every non-closed
subset~$A$ of~$X$, there exists a sequence in~$A$ converging to some
point~$x$ in~$\closure A\setminus A$. Contrast this with the situation for
a Fréchet space: for every non-closed subset~$A$ and \emph{every}
point~$x$ in~$\closure A\setminus A$, there exists a sequence in~$A$
converging to~$x$. For a Fréchet space, one can choose the point~$x$ in
advance, whereas for a sequential space one cannot. Thus every Fréchet
space is sequential, but not conversely (see
Example~\ref{eg:sequential-frechet}).
\end{rem}

We can say right now that none of our classes are invariant under inverse 
images of proper mappings. Indeed, $\beta\nats$ is not sequential and yet 
$f\maps\beta\nats\to\set{0}\maps x\mapsto 0$ is a 
\index{proper mapping}proper mapping (i.e., it is onto, continuous, 
closed and has compact fibres). Likewise no class is closed with respect 
to images under continuous mappings: for if $X$~denotes the set 
$\beta\nats$ with the discrete topology, then $X$~is metrisable and the 
identity mapping $\iota\maps X\to\beta\nats$ is continuous.

%\medskip

We have defined some classes of spaces in which the topology (in 
particular, closed sets) may be described by sequences. We now generalise 
a little and look at spaces in which the closed sets may be described by 
countable sets, in the following sense.

\begin{defn}[countable tightness]\index{tightness|see{countable tightness}}
\label{defn:counttight}
Let $X$~be a topological space. We say that $X$~has 
\define{countable tightness} if whenever $A$~is a subset of~$X$ that 
contains the closure of all its countable subsets, then $A$~is closed.
\end{defn}

The definition of countable tightness that we have given is analogous to
that of sequential space, in the sense that if you replace `countable
subsets' by `sequences' and `closure' by `limits' in the present
definition, then you arrive at a definition of sequential space. Thus
every sequential space, and hence every first countable space, has
countable tightness.

There is an alternative characterisation of tightness, which is analogous
in the same way to the definition of a Fréchet space. When the word
\emph{tightness} was first introduced (by \citet{Arhangelskij68}), the
class of spaces that we just defined would have been said to have
countable \emph{weak tightness}. The next result is a special case of the
fact that tightness and weak tightness are one and the same thing. In the
light of Remark~\ref{rem:frechetsequential}, this is somewhat surprising.

\begin{lem}\label{lem:counttight-equiv}
Let $X$~be a topological space. Then $X$~has countable tightness if
and only if whenever a point~$x\in X$ belongs to the closure of
some~$A\sub X$, then there exists a countable subset~$C$ of\/~$A$
such that $x$~belongs to~$\closure C$.
\end{lem}

\begin{proof}
  Suppose that $X$~is a space such that for all $x\in\closure A$ there
  exists countable $C\sub A$ such that $x\in\closure C$. This implies
  that if $A$~is a subset of~$X$ that contains the closure of all its
  countable subsets, then $A$~contains all the points of its closure,
  so $A=\closure A$ and $A$~is closed. Therefore $X$~has countable
  tightness.

  Now suppose conversely that $X$~has countable tightness, and let
  $A$~be a subset of~$X$. Let
  \begin{displaymath}
    E\defeq\Union\set{\closure
      C\st\text{$C$ countable, $C\sub A$}}\text.
  \end{displaymath}
  We want to show that $\closure A\sub E$, because then every point in
  $\closure A$ belongs to~$E$, so it lies in the closure of some
  countable subset of~$A$, as desired. Clearly $A\sub E$, so to show
  that $\closure A\sub E$ it is sufficient to show that $E$~is closed:
  to do this we shall use countable tightness. Let $D$~be a countable
  subset of~$E$, and let $x\in\closure D$. Observe that for every
  $y\in D$ we have $y\in\closure{C_y}$ for some countable
  subset~$C_y$ of~$A$. Let $C\defeq\Union_{ y\in D}C_y$; since
  $D$~is countable, $C$~is countable, and it is certainly a subset
  of~$A$. Now let $U$~be any open neighbourhood of~$x$. Then from
  the fact that $x\in\closure D$, there exists some $y$ in
  $D\intersect U$. Since $U$~is open it is a neighbourhood of~$y$,
  and since $y\in\closure{C_y}$ it follows that $U$~meets~$C_y$,
  hence $U$~meets~$C$. This shows that $x\in\closure C$, which
  implies that $x\in E$. Therefore $\closure D\sub E$, thus
  $E$~contains the closure of all its countable subsets. Since
  $X$~has countable tightness, it follows that $E$~is closed. Now we
  see that~$E$ is closed and $A\sub E$, which implies that $\closure
  A\sub E$. It follows that, for any $x\in\closure A$, there exists
  some countable subset~$C$ of~$A$ such that $x\in\closure C$.
\end{proof}

Lemma~\ref{lem:counttight-equiv} was stated in a footnote by 
\citet{Arhangelskij70}. This proof is my answer to an exercise 
in~\citet[Prob.~1.7.13]{Engelking89}.

It is easy to show that every
\index{hereditarily separable space}%
hereditarily separable space~$X$ has countable tightness: for any 
$A\sub X$, let $C$~be a countable dense subset. This is dense in the sense 
that $\closure[A]C=A$, i.e., $\closure C\intersect A=A$, whence 
$A\sub\closure C$. Now if $A$~contains the closure of all its countable 
subsets then $\closure C\sub A$, it follows that $A=\closure C$ and $A$~is 
closed. Thus a hereditarily separable space has countable tightness. In 
particular, every countable space has countable tightness.

%\medskip

Another class of spaces can also be defined in terms of a `countable' 
property of closed sets, as follows.

\begin{defn}[perfect space]
Let $X$~be a topological space, and let $A$~be a subset of~$X$. We say 
that $A$~is a \define[gdelta@$G_\delta$~set]{$G_\delta$~set} if it can be 
written as the intersection of a countable family of open sets. We say 
that $A$~is an \define[fsigma@$F_\sigma$~set]{$F_\sigma$~set} if it can be 
written as the union of a countable family of closed sets.

The topological space~$X$ is said to be \define[perfect space]{perfect} if 
every closed subset of~$X$ is a~$G_\delta$. A space that is both perfect 
and normal is called \define{perfectly normal}.
\end{defn}

It is clear that a space is perfect iff every open set is an~$F_\sigma$. 
We discuss perfect spaces in Section~\ref{chap:perfect}.

%%%%%%%%%%%%%%%%%%%%%%%%%%%%%%%%%%%%%%%%%%%%%%%%%%%%%%%%%%%%%%%%%%%%%
%%%%%%%%%%%%%%%%%%%%%%%%%%%%%%%%%%%%%%%%%%%%%%%%%%%%%%%%%%%%%%%%%%%%%
%%%%                                                             %%%%
%%%%                      CHAPTER TWO                            %%%%
%%%%                                                             %%%%
%%%%%%%%%%%%%%%%%%%%%%%%%%%%%%%%%%%%%%%%%%%%%%%%%%%%%%%%%%%%%%%%%%%%%
%%%%%%%%%%%%%%%%%%%%%%%%%%%%%%%%%%%%%%%%%%%%%%%%%%%%%%%%%%%%%%%%%%%%%

\section{First Countable Spaces}\label{chap:firstcount}
\index{first axiom of countability|(}

\subsection*{The Basics}

This section is another introductory one, in which we
review the properties of first countable spaces. This in
preparation for the real work of looking at sequential spaces in
Section~\ref{chap:sequential}. There are two main points of interest.
The first is that, in first countable spaces, the Hausdorff
separation axiom is equivalent to the uniqueness of sequential
limits (Proposition~\ref{prop:firstcount-hausdorff}); this result is no
longer true for the class of Fréchet spaces. The second point to
note is that in passing from metric spaces to first countable
spaces we lose the equivalence of compactness and sequential
compactness. We give an example of a first countable sequentially
compact space that is not compact in
Example~\ref{eg:firstcount-compact}. Finally, for comparison with the
main theorems of
Sections~\ref{chap:sequential}--\ref{chap:frechet}, we shall state
Ponomarev's characterisation of the $T_0$~first countable spaces
in Theorem~\ref{thm:firstcount-characterisation}.

Naturally one asks what are the operations on first countable
spaces that yield a first countable space. The operation of taking
subspaces is one of them, according to the following result. We
omit the (easy) proof.
\begin{lem}\label{lem:firstcount-hereditary}
  \index{subspace!of first countable space}
  Any subspace of a first countable space is again first
  countable.\qed
\end{lem}

Not every quotient of a first countable space is first countable;
this will follow from results in Section~\ref{chap:sequential}.

We shall now give an example to show that the product of first
countable spaces need not be first countable.
\begin{example}
  \index{product!of first countable space}%
  \label{eg:firstcount-productive}%
  Consider the space~$\real^\real$, i.e., the Tychonoff product of
  $\continuum$~copies of the real line. Since $\real$ is a metric
  space it is first countable; we shall show that $\real^\real$ is
  \emph{not} first countable. Suppose for a contradiction that
  $\Bee=\set{U^{(i)}\st i\in\nats}$~is a countable base of open
  neighbourhoods of~$0$ in~$\real^\real$. For each natural
  number~$i$, since $U^{(i)}$~is open there exists a basic open
  set~$\prod_{r\in\real}W^{(i)}_r$ such
  that~$0\in\prod_{r\in\real}W^{(i)}_r\sub U^{(i)}$.
  Since~$\prod_{r\in\real}W^{(i)}_r$ is basic and open, there exists a
  finite subset~$R^{(i)}$ of~$\real$ such that~$W^{(i)}_r=\real$ for
  every~$r\in\real\setminus R^{(i)}$. Now $R\defeq\Union\set{R^{(i)}\st
    i\in\nats}$ is a countable, hence proper, subset of~$\real$.
  Therefore there exists~$r\in\real$ with~$r\notin R$, and this
  $r$~satisfies $W^{(i)}_r=\real$ for every~$i\in\nats$. In
  particular~$W^{(i)}_r\nsubseteq(-1,1)$, for each~$i\in\nats$. Hence
  \begin{equation}\label{eq:1}
    \prod_{s\in\real}W^{(i)}_s\nsubseteq(-1,1)\setprod\prod_{s\neq
      r}\real=\pi_r{}^{-1}(-1,1)\text,
  \end{equation}
  from which it follows that $U^{(i)}\nsubseteq\pi_r{}^{-1}(-1,1)$, for
  every~$U^{(i)}\in\Bee$. But $\pi_r{}^{-1}(-1,1)$ is an open
  neighbourhood of~$0$ in~$\real^\real$, which contradicts the fact
  that $\Bee$~is a neighbourhood base at~$0$. This shows that
  $\real^\real$ is not first countable.
\qed
\end{example}

The argument above is adapted from \citet[Thm.~3.6]{Kelley55}.
Note that in this example we also see that an uncountable product
of metric spaces need not be metrisable. However, it is easy to
show that a \emph{countable} product of metric spaces is
metrisable. That result carries across to first countable spaces
in the sense that a countable product of first countable spaces is
first countable, as we now show.
\begin{prop}\index{product!of first countable space}
  \label{prop:firstcount-countableproduct}
  Suppose that $X_n$ is a first countable space for each~$n$ in~$\Z$,
  and let $X$~denote the Tychonoff product $\prod_{n\in\Z}X_n$. Then
  $X$~is a first countable space.
\end{prop}
\begin{proof}
  Let $x$~belong to~$X$. For each $n\in\Z$, let
  $\set{U^{(i)}_n\st i\in\nats}$ be a countable base for the neighbourhood
  system at~$x_n$. Let
  \begin{gather*}
    N^{(i)}_n\defeq U^{(i)}_n\setprod\prod_{m\neq n}X_m
    =\pi_n{}^{-1}(U^{(i)}_n)\text,\qquad\text{for $i\in\nats$, $n\in\Z$,}\\
    {\Bee\defeq\set{\Intersection_{n\in F}N^{(i(n))}_n\st
        \text{$i(n)\in\nats$, $F\sub\Z$ is finite}}}\text.
  \end{gather*}
  Then $\Bee$~is a countable family of open subsets of~$X$. Now, let
  $N$~be a basic open neighbourhood of~$x$. Then there exists a
  finite subset~$F$ of~$\Z$ and a family $\set{V_n\st n\in\Z}$ of
  sets such that $V_n$~is open in~$X_n$ for all~$n$, with $V_n=X_n$
  for all $n\in\Z\setminus F$ and $N=\prod_{n\in\Z}V_n$. Pick any~$n$
  in~$F$; then $V_n$~is a neighbourhood of~$x_n$ in~$X_n$. So for
  each $n\in F$ there exists some natural number~$i(n)$ such that
  $x_n\in U^{(i(n))}_n\sub V_n$. Therefore $x\in\Intersection_{n\in
    F}N^{(i(n))}_n\sub N$. This shows that $\Bee$~is a base for the
  neighbourhood system at~$x$. Hence $X$~is first countable.
\end{proof}

\subsection*{For First Countable Spaces Only}

Now we present a result that is \emph{not} available for Fréchet
spaces (cf.\ Example~\ref{eg:frechet-hausdorff}). This
characterisation of the Hausdorff property in the realm of first
countable spaces is part of topological folklore (see for instance
\citet[Chap.~III, Ex.~3A]{Hu64}).
\begin{prop}\label{prop:firstcount-hausdorff}\index{Hausdorff property}
  Let~$X$ be a first countable space. Then $X$~is Hausdorff if and
  only if sequential limits in~$X$ are unique, that is whenever
  $\seq(x_i)$ converges to~$x$ and to~$y$, then~$x= y$.
\end{prop}
\begin{proof}
  Suppose first that $X$~is a Hausdorff space. Let~$\seq(x_i)$ be a
  sequence converging to some point~$x$, and let $y\neq x$. Then
  there exist disjoint open sets $U$,~$V$ such that $x\in U$ and $y\in
  V$. Since $x_i\to x$ there exists a natural number $i_0$ such that
  $x_i\in U$ whenever $i\gee i_0$. Since $U$,~$V$ are disjoint it
  follows that $x_i\notin V$ whenever $i\gee i_0$. Therefore
  $\seq(x_i)$ cannot possibly converge to~$y$. This shows that
  sequential limits in~$X$ are unique.

  Conversely, suppose that $X$~is not a Hausdorff space. Then there
  exist two distinct points $x$,~$y$ in~$X$ such that every
  neighbourhood of~$x$ meets every neighbourhood of~$y$.
  Let~$\Bee=\set{U_{i}\st i\in\nats}$ be a base of open neighbourhoods
  at~$x$, and let~$\set{V_{i}\st i\in\nats}$ be a base of open
  neighbourhoods at~$y$. Then for each~$i\in\nats$ we have
  $U_{i}\intersect V_{i}\neq\emptyset$, so there exists some
  point~$x_i$ in~$U_{i}\intersect V_{i}$. Since~$\Bee$ is a base
  at~$x$, the sequence~$\seq(x_i)$ converges to~$x$. Similarly
  $\seq(x_i)$ converges to~$y$. But by assumption~$x\neq y$, thus
  sequential limits in~$X$ are not unique.
\end{proof}
We can see that the topology of first countable spaces is very
largely determined by their convergent sequences. By
Lemma~\ref{lem:firstcount-frechet} the closure operation can be
described by sequences, and every first countable space is
sequential, so the open sets are determined by sequences as well.
In Lemma~\ref{lem:sequential-conts} we shall see that continuity too
can be characterised by sequences, and
Proposition~\ref{prop:firstcount-hausdorff} just proved shows that the
same can be said for the Hausdorff property. Thus, most of the
properties of first countable spaces are determined by knowledge
of their convergent sequences. However, there is a substantial fly
in the ointment: compactness.

\subsection*{An Ordinal Space}

For sequential Hausdorff spaces, sequential compactness is
equivalent to countable compactness (see
Proposition~\ref{prop:sequential-compact}). In particular, a first
countable compact Hausdorff space is sequentially compact.
However, this result has no converse: a first countable
sequentially compact Hausdorff space need not be compact. The
following example is standard, it appears in
\citet[Sec.~3.7]{Munkres75}, for instance.
\begin{example}\index{omega@$\omega_1$}%
  \label{eg:firstcount-compact}\index{compact space|(}
  \index{sequentially compact space|(}
  Let $\omega_1$~be the least uncountable ordinal number, and define
  \begin{displaymath}
    W_0\defeq\set{x\st\text{$x$~is an ordinal and $x<\omega_1$}}\text.
  \end{displaymath}
  Then $W_0$~is a well-ordered set with ordinality~$\omega_1$; the
  topology on~$W_0$ is the order topology. We are going to show that
  $W_0$~is a sequentially compact Hausdorff first countable space that
  is not compact. To show that $W_0$~is not compact we show that it
  is a dense proper subspace of a compact Hausdorff space. To this
  end, let $W\defeq W_0\union\set{\omega_1}$, and let $W$~have the
  order topology. First we establish some notation. Let
  $-\infty\defeq\min W$, and given $x\in W$ write
  $[-\infty,x)$,~$(x,\omega_1]$ for the sets $\set{y\in W\st y<x}$,
  $\set{y\in W\st y>x}$ respectively.

  By definition these sets form a subbase for the order topology
  on~$W$. When $x<y$ we define $(x,y)$,~$[x,y)$, $(x,y]$ and~$[x,y]$
  in the obvious way. Note that whenever $x<y$ we have
  $(x,y]=[-\infty,y+1)\intersect(x,\omega_1]$, so that $(x,y]$ is
  open. Now, every neighbourhood~$N$ of $\omega_1$ contains a set of
  the form $(x,\omega_1]$ for some $x\in W_0$, hence $N$~meets~$W_0$.
  Therefore $\omega_1\in\closure{W_0}$, so that $W_0$~is dense in~$W$.

  To show that $W$~is compact, let $\Eee$~be an open cover for~$W$ and
  consider the set
  \begin{displaymath}
    S\defeq\set{x\in W\st\vphantom{_j}\text{$[-\infty,x]$ has a finite
        subcover in~$\Eee$}}\text.
  \end{displaymath}
  Suppose for a contradiction that $S\neq W$. Then $W\setminus S$ is
  non-empty, so it has a least element, $x$~say. Since~$\Eee$ is a
  cover, there exists~$U_0$ in~$\Eee$ with $x\in U_0$. Clearly
  $-\infty\in S$, so $x\neq-\infty$. Since $U_0$~is open we have
  $x\in[-\infty,y)\intersect(z,\omega_1]\sub U_0$ for some $y,z\in W$,
  whence $(z,x]\sub U_0$. Since $z<x$ we have $z\in S$, so
  $[-\infty,z]$ can be covered by a finite subcollection
  $\set{U_1,U_2,\dots,U_n}$ of~$\Eee$. But now
  $\set{U_0,U_1,\dots,U_n}$ is a finite subcollection of~$\Eee$ that
  covers $[-\infty,x]$, which contradicts the fact that $x\notin S$.
  Therefore $S=W$ after all, i.e., $W$~has a finite subcover
  in~$\Eee$. This shows that $W$~is compact.

  Now we show that $W$ is Hausdorff: if $x$,~$y$ are distinct points
  of~$W$, then we may assume without loss of generality that $x<y$.
  Then $U\defeq[-\infty,x+1)$ and $V\defeq(x,\omega_1]$ are disjoint
  open sets with $x\in U$ and $y\in V$, whence $W$~is a Hausdorff
  space. Hence $W_0$~is a non-closed subspace of the compact
  Hausdorff space~$W$, so $W_0$~is not compact. As a subspace of~$W$,
  \ $W_0$~is also a Hausdorff space.

  Next we can show that $W_0$~is first countable. Pick~$x\in W_0$.
  If $x=-\infty$ then $\set{x}$ is a countable neighbourhood base.
  Otherwise, the collection $\Bee\defeq\set{(y,x]\st y<x}$ is
  countable because $\Bee$ is well-ordered by inclusion, with
  ordinality $x<\omega_1$. Since $\Bee$~is a base for the
  neighbourhood system at~$x$, this shows that $W_0$~is first
  countable.

  We want to show that $W_0$~is sequentially compact. First we note
  the following fact: {\itshape If $C$ is a countable subset
    of\/ $W_0$, then $C$ is bounded above.} For, let $x$~be an
  element of~$C$. Then since $x<\omega_1$, it follows that
  $[-\infty,x]$ is countable. Therefore the set $B\defeq\Union_{x\in
    C}[-\infty,x]$ is a countable union of countable sets, hence
  countable. It follows that $B\neq W_0$, since $W_0$ is uncountable.
  Therefore there exists some $y\in W_0\setminus B$. Since $y\notin
  B$ we must have $y>x$ for all~$x\in C$, that is, $y$~is an upper
  bound for~$C$.

  To show that $W_0$~is sequentially compact, let $\seq(x_i)$~be a
  sequence in~$W_0$. Consider the set
  \begin{displaymath}
    F\defeq\set{\smash{x_i\st\text{$x_i\gee x_j$
          whenever $j\gee i$}}}\text.
  \end{displaymath}
  If $F=\emptyset$, then let $I\defeq\set{1}$. Otherwise, let
  $x\defeq\min F$ and let $I\defeq\set{i\in\nats\st x_i=x}$. In the
  case that $I$~is infinite, we can arrange~$I$ into an increasing
  sequence $\seq({i(r)}\st r\in\nats)$, and then we see that
  $\seq(x_{i(r)})$ is a constant, hence convergent, subsequence of
  $\seq(x_i)$.

  It only remains to consider the case that $I$~is finite. Then let
  $i(1)\defeq\max I+1$, and consider~$x_{i(1)}$. Then $x_{i(1)}\notin
  F$, so there exists some natural number $i(2)>i(1)$ with
  $x_{i(2)}>x_{i(1)}$. Continuing in this way we can pick out a
  strictly increasing subsequence $\seq(x_{i(r)})$ of $\seq(x_i)$.
  Now the set $C=\set{x_{i(r)}\st r\in\nats}$ is countable, hence
  bounded above. Thus the set of all upper bounds of~$C$ is
  non-empty, let $y$~be its least element. Let $N$~be a neighbourhood
  of~$y$; then $N$~contains~$(z,y]$ for some $z\in W_0$. Therefore
  $z<y$, so $z$~is not an upper bound for~$C$, whence $z<x_{i(s)}\lee
  y$ for some natural number~$s$. Since the sequence $\seq(x_{i(r)})$
  is increasing it follows that $z<x_{i(r)}\lee y$, hence $x_{i(r)}\in
  N$, whenever $r\gee s$. This shows that $x_{i(r)}\to y$, thus
  $\seq(x_i)$ has a convergent subsequence. Therefore $W_0$~is
  sequentially compact.%
  \index{sequentially compact space|)}\index{compact space|)}
\qed
\end{example}

We finish off this section by stating a theorem that characterises
the first countable spaces. Recall that a continuous mapping
$f\maps X\to Y$ is said to be \define[open mapping]{open} if it
maps open sets in~$X$ to open sets in~$Y$.
\begin{thm}[\citet{Ponomarev60}]
  \label{thm:firstcount-characterisation}
  If $f$~is an open mapping from a metric space~$X$ onto a
  topological space~$Y$, then $Y$~is first countable. Conversely, if\/~$Y$
  is a first countable space then for some metric space~$X$
  there exists an open mapping~$f$ from~$X$ onto~$Y$.
\end{thm}
This theorem will not be proved here. It is included mainly for
comparison with our other characterisation theorems
(Theorems~\ref{thm:sequential-characterisation}
and~\ref{thm:frechet-characterisation}).%
\index{first axiom of countability|)}

%%%%%%%%%%%%%%%%%%%%%%%%%%%%%%%%%%%%%%%%%%%%%%%%%%%%%%%%%%%%%%%%%%%%%
%%%%%%%%%%%%%%%%%%%%%%%%%%%%%%%%%%%%%%%%%%%%%%%%%%%%%%%%%%%%%%%%%%%%%
%%%%                                                             %%%%
%%%%                      CHAPTER THREE                          %%%%
%%%%                                                             %%%%
%%%%%%%%%%%%%%%%%%%%%%%%%%%%%%%%%%%%%%%%%%%%%%%%%%%%%%%%%%%%%%%%%%%%%
%%%%%%%%%%%%%%%%%%%%%%%%%%%%%%%%%%%%%%%%%%%%%%%%%%%%%%%%%%%%%%%%%%%%%

\section{Sequential Spaces}\label{chap:sequential}
\index{sequential space|(}

\subsection*{Sequences in Sequential Spaces}

At this point we begin our work in earnest by examining
the sequential spaces (recall Definition~\ref{defn:sequential}). Our
ultimate goal in this section is
Theorem~\ref{thm:sequential-characterisation}, which states that every
sequential space is the quotient of a metric space, and \emph{vice
versa}. This theorem is quite remarkable. First, quotient spaces
of metric spaces are easy to form so we have a plentiful supply of
sequential spaces. Indeed this shows that many of the most natural
examples of topological spaces are sequential. Second, the theorem
gives a characterisation of the class of sequential spaces in
terms that do not involve sequential convergence.

How well do sequences describe the topology in a sequential space?
Open sets and closed sets admit simple descriptions, by
definition. The closure operation is not so well-behaved; we shall
give an example of a sequential space that is not Fréchet in
Example~\ref{eg:sequential-frechet}. However, the continuity of a
function \emph{can} be characterised in terms of sequences, as
follows.\par
\begin{lem}\index{continuous mapping!of sequential space}
  \label{lem:sequential-conts}
  Let $X$~be a sequential space, let $Y$~be an arbitrary topological
  space and let $f\maps X\to Y$. Then $f$~is continuous if and only
  if it preserves limits of sequences, i.e., whenever $\seq(x_i)$ is a
  sequence converging in~$X$ to~$x$, then $\seq({f(x_i)})$ converges
  in~$Y$ to~$f(x)$.
\end{lem}
\begin{proof}
  Suppose that $f$~is continuous, and let $\seq(x_i)$ be a sequence
  in~$X$ with $x_i\to x$. Let $N$~be a neighbourhood of~$f(x)$.
  Since $f$~is continuous, $f^{-1}(N)$ is a neighbourhood of~$x$.
  Therefore there exists some natural number~$i_0$ such that $x_i\in
  f^{-1}(N)$ whenever $i\gee i_0$. This implies that $f(x_i)\in N$
  whenever $i\gee i_0$, therefore $f(x_i)\to f(x)$.

  Now suppose conversely that $f(x_i)\to f(x)$ whenever $x_i\to x$.
  Let $U$~be an open subset of~$Y$. Let $x\in f^{-1}(U)$ and suppose
  that $x_i\to x$. Then $f(x_i)\to f(x)$, where $f(x)\in U$. Now
  $U$~is open, so there exists a natural number~$i_0$ such that
  $f(x_i)\in U$ whenever $i\gee i_0$. Hence $x_i\in f^{-1}(U)$
  whenever $i\gee i_0$. This shows that $f^{-1}(U)$ is sequentially
  open. Since $X$~is sequential, it follows that $f^{-1}(U)$ is open,
  thus $f$~is continuous.
\end{proof}

In general compactness neither implies nor is implied by sequential
compactness. We gave an example (Example~\ref{eg:firstcount-compact})
showing that, even for Hausdorff first countable spaces, sequential
compactness does not imply compactness. The space $\beta\nats$~is an
example of a compact Hausdorff space that is not sequentially
compact. However, it \emph{is} true that a compact Hausdorff
sequential space is sequentially compact. In order to prove this we
shall actually show that, in Hausdorff sequential spaces, sequential
compactness is equivalent to countable compactness.  This result was
first noted by \citet[Prop.~1.10]{Franklin65}. Our proof is based on
the sketch in \citet[Thm.~3.10.31]{Engelking89}
\begin{prop}\index{compact space}\index{countably compact space}
  \index{sequentially compact space}
  \label{prop:sequential-compact}
  Let $X$~be a Hausdorff sequential space. Then $X$~is sequentially
  compact if and only if it is countably compact.\footnote{I am
    grateful to Alexander Gouberman for pointing out an error in an
    earlier version of this proposition.}
\end{prop}
\begin{proof}
  Sequential compactness always implies countable compactness: use the
  obvious proof by contradiction. To prove the converse, let $X$~be a
  countably compact Hausdorff sequential space. Let $\seq(x_i)$ be a
  sequence in~$X$. Let $A\defeq\set{x_i\st i\in\nats}$. Since $X$~is
  countably compact, the sequence $\seq(x_i)$ has a cluster point~$x$
  in~$X$---see e.g., \citet[IV.5.D)]{Nagata85}. Thus, for every
  neighbourhood~$N$ of~$x$, and for all natural numbers~$j$, there
  exists some~$i\gee j$ such that $x_i$~belongs to~$N$. If $x_i=x$ for
  infinitely many~$i$, then we have a convergent subsequence. If not,
  then $x$~must be an element of~$\closure{(A\setminus\set{x})}$.
  Therefore ${A\setminus\set{x}}$ is not closed. Since $X$~is a
  sequential space, it follows that ${A\setminus\set{x}}$~is not
  sequentially closed. Hence there exists a sequence $\seq(y_r)$ in
  ${A\setminus\set{x}}$ that converges to some point~$y$ with $y\notin
  {A\setminus\set{x}}$. By virtue of the Hausdorff property, we may
  assume that all of the points $y_r$ are distinct. Now $y_r\in A$ for
  each natural number~$r$, so there exists some $j(r)$ such that
  $y_r=x_{ j(r)}$. Now we delete those $x_{ j(r)}$ that occur `too
  late', by defining a function $i\maps\nats\to\nats$ in the following
  inductive manner:
  \begin{displaymath}
    i(1)\defeq j(1)\text,\qquad
    i(r+1)\defeq j\bigl(\min\set{s>r\st\smash{j(s)>i(r)}}\bigr)\text.
  \end{displaymath}
  This is contrived so that $i$~is monotone increasing, thus
  $\seq({x_{i(r)}})$ is a subsequence of $\seq(x_i)$. Also, it is
  subsequence of $\seq(y_r)$, so it converges to~$y$. This shows that
  $X$~is sequentially compact.
\end{proof}

\subsection*{Quotients of Sequential Spaces}

Recall that a continuous onto mapping $f\maps X\to Y$ is called a
\define{quotient mapping} if $V$~is open in~$Y$ whenever
$f^{-1}(V)$~is open in~$X$. For example, the projection mapping
from a product space onto one of its coordinate spaces is a
quotient mapping, a fact that we shall use in
Example~\ref{eg:sequential-frechet}.
\begin{lem}[{\citet[Prop.~1.2]{Franklin65}}]
  \label{lem:sequential-quotient}
  \index{quotient mapping!of sequential space}%
  Let $X$~be a sequential space, and let $f$~be a quotient mapping
  from~$X$ onto a topological space~$Y$. Then $Y$~is a sequential
  space.\qed
\end{lem}
This is proved by showing that if~$U$ is sequentially open in~$Y$,
then $f^{-1}(U)$ is sequentially open in~$X$.
\begin{cor}\index{proper mapping!of sequential space}
  \label{cor:sequential-properimage}
  The class of sequential spaces is closed with respect to taking
  images under proper mappings.\qed
\end{cor}

With the help of Lemma~\ref{lem:sequential-quotient} we can give an
example of a sequential space that is not Fréchet. This shows
that, in sequential spaces, sequences do not in general describe
the closure operator. 

%\pagebreak[2]
\begin{example}[{\citet[Eg.~2.3]{Franklin65}}]
  \label{eg:sequential-frechet}\index{Fréchet space}
  Let $L\defeq\real\setminus\set{0}$ be the set of real numbers with
  $0$~removed, and let $M\defeq\set{0}\union\set{1/i\st i\in\nats}$.
  Let $Y\defeq(L\setprod\set{0})\union(M\setprod\set{1})$ have its
  usual (metric) topology as a subspace of $\real^2$. In particular,
  $Y$~is a sequential space. Now let $f\maps Y\to X$ be the
  projection of~$Y$ onto its first coordinate space, the
  set~$X\defeq\set{0}\union L$. Observe that $X$~is the set of all
  real numbers. The topology on~$X$ is generated by the usual
  topology plus all sets of the form $\set{0}\union U$ where $U$~is
  open in~$\real$ and contains $\set{1/i\st i\in\nats}$. As a
  quotient of a sequential space, $X$~is sequential by
  Lemma~\ref{lem:sequential-quotient}.

  Now we shall show that $X$~is not a Fréchet space. Let $A\defeq
  X\setminus M$. We show that $0\in\closure A$ but that no sequence
  in~$A$ converges to~$0$. Let $N$~be a neighbourhood of~$0$. Then
  $N$~contains a point $1/i$ for some $i\in\nats$. Away from~$0$, the
  topology of~$X$ is locally `as usual'. At any rate, there exists
  some $\delta>0$ such that $(1/i-\delta,1/i+\delta)\sub N$. Now
  there exists some irrational~$x$ with $\abs{1/i-x}<\delta$. Then
  $x\in N$ and, $x$~being irrational, $x\in A$. Hence $A\intersect
  N\neq\emptyset$. As this was true for any neighbourhood of~$0$,
  this shows that $0\in\closure A$. Now suppose that $\seq(x_i)$ is
  any sequence in~$A$. For each $i\in\nats$, we have $x_i\neq 0$, and
  therefore $\inf_{ j\in\nats}\abs{x_i-1/ j}>0$. Consider the
  sequence $\seq(\delta_i)$, where
  $\delta_i\defeq\inf_{ j\in\nats}\abs{x_i-1/ j}$, and let
  \begin{displaymath}
    N\defeq\set{0}\union\Union_{i\in\nats}(1/i-\delta_i,1/i+\delta_i)
    \text.
  \end{displaymath}
  Then $N$~is a neighbourhood of~$0$ that contains none of the members
  of the sequence $\seq(x_i)$. Therefore $\seq(x_i)$ cannot possibly
  converge to~$0$. This shows that no sequence in~$A$ converges
  to~$0$, so $X$~is not a Fréchet space.
\qed
\end{example}

We now look at how well the property of sequentiality behaves
under the operations that we like to use to in forming new
topological spaces. Granting that quotients of sequential spaces
are sequential, the class of sequential spaces does not in fact
possess many other nice closure properties. With
Example~\ref{eg:sequential-productive} we shall show that the product
of two sequential spaces need not be sequential. First, we are
going to give an example of a sequential space with a
non-sequential subspace.
\begin{example}[{\citet[Eg.~1.8]{Franklin65}}]
  \label{eg:sequential-hereditary}\index{subspace!of sequential space}
  Let $X$~be the space in Example~\ref{eg:sequential-frechet}, and let
  $A$~be as above. Let $Z\defeq A\union\set{0}$. We show that $Z$~is
  a subspace of~$X$ that is not sequential. Recall that no sequence
  in~$A$ converges to~$0$. Therefore the only sequences in~$Z$ that
  converge to~$0$ are those that are eventually constant. Hence the
  set~$\set{0}$ is sequentially open in~$Z$. However, $\set{0}$~is
  not open in~$Z$, for $f^{-1}\set{0}=\set{0}\setprod\set{1}$ is not
  open in~$Y$, as any neighbourhood of $\set{0}\setprod\set{1}$
  contains some point $(1/i,1)$. Thus $Z$~is not a sequential space.
\qed
\end{example}

Although a general subspace of a sequential space might not be
sequential, a closed subspace must be.
\begin{lem}[{\citet[Prop.~1.9]{Franklin65}}]
  \index{subspace!of sequential space}
  \label{lem:sequential-closedsubspace}
  A closed subspace of a sequential space is again sequential.\qed
\end{lem}

It should be noted that, in Example~\ref{eg:sequential-hereditary},
the space $X$~is a sequential space that is not Fréchet. This is
not a coincidence; in fact any example of a non-hereditarily
sequential space must be of this type:
\begin{lem}[{\citet[Prop.~7.2]{Franklin67}}]
  \label{lem:hereditarilysequential-frechet}
  \index{subspace!of sequential space}%
  Let $X$~be a sequential space. If every subspace of\/~$X$ is
  sequential, then $X$~is Fréchet.\qed
\end{lem}
Recalling Remark~\ref{rem:frechetsequential}, the proof is easy: for
a non-closed subset $A$, take a point $x\in\closure A\setminus A$,
then consider the subspace $A\union\set{x}$.

We conclude our short investigation into the closure properties of
the class of sequential spaces with the following example. This
nice argument is due to T. K. Boehme, cited in \citet{Franklin66}.
\begin{example}\index{product!of sequential space}
  \label{eg:sequential-productive}%
  Let $X$~be any non-Hausdorff Fréchet space with unique sequential
  limits. (By Example~\ref{eg:frechet-hausdorff} below, such an~$X$
  exists; at present we omit the construction.) We show that the space
  $X\setprod X$ is not sequential. Since $X$~is not Hausdorff, the
  diagonal $\Delta\defeq\set{(x,x)\st x\in X}$ is not closed
  in~$X\setprod X$. On the other hand, $\Delta$~is sequentially
  closed, because sequential limits are unique. Hence $X\setprod
  X$~is not a sequential space.
\qed
\end{example}

\subsection*{The Main Theorem}

We are now going to state and prove our characterisation theorem,
that identifies every sequential space as the image of a metric
space under a quotient mapping.
\begin{thm}[{\citet{Franklin65}}]
  \label{thm:sequential-characterisation}\index{quotient mapping}
  If $f$~is a quotient mapping from a metric space~$X$ onto a
  topological space~$Y$, then $Y$~is sequential. Conversely, if\/
  $Y$~is a sequential space, then for some metric space~$X$ there
  exists a quotient mapping~$f$ from~$X$ onto~$Y$.
\end{thm}
The first half of Theorem~\ref{thm:sequential-characterisation} is
immediate from Lemma~\ref{lem:sequential-quotient}. To prove the
second half we must find a suitable metric space~$X$. We
construct~$X$ as a topological sum of convergent sequences.

\begin{rem}\label{rem:convergent-sequence}
  By a \define{convergent sequence}, we mean a sequence $\seq(x_i)$ of
  points together with a single limit~$x$ (in this we are following
  \citet{Franklin65}). We consider the set $X\defeq\set{x_i\st
    i\in\nats}\union\set{x}$ in its natural topology: each of the
  points~$x_i$ is isolated, and the neighbourhoods of~$x$ are
  precisely the complements of the finite subsets of~$X$. Topologised
  in this way, $X$~is the Alexandroff one-point compactification of
  the discrete space $\set{x_i\st i\in\nats}\setminus\set{x}$. Since
  the family of finite subsets of a countable set is countable, we see
  that $X$~is a second countable compact Hausdorff space, so by
  Urysohn's metrisation theorem, $X$~is metrisable.
\end{rem}

\begin{proof}[Proof of Theorem~\ref{thm:sequential-characterisation}]
  Let $Y$~be a sequential space, and let $\Ess$~be the set of all
  convergent sequences in~$Y$. For each sequence $S\maps\nats\to
  Y\maps i\mapsto x_i$ in~$\Ess$, let $L_S\defeq\set{x\in Y\st x_i\to
    x}$. For each $x\in L_S$, consider the space
  $S(\nats)\union\set{x}$ endowed with the metrisable topology
  described above. We wish to define $X$~to be the disjoint
  topological sum of all such spaces. However, the spaces
  $S(\nats)\union\set{x}$ as we have described them are not, \emph{de
    facto}, disjoint. Before forming the sum we need disjoint
  homeomorphic copies of these spaces. This can be achieved
  (cf.~\citet[p.~75]{Engelking89}) by taking the product with a
  one-point discrete space (twice): let
  $X_{(S,x)}\defeq(S(\nats)\union\set{x})\setprod\set{S}\setprod\set{x}$.
  Then the spaces $X_{(S,x)}$ are disjoint, so we may form their
  topological sum:
\begin{displaymath}
X \defeq 
%\Dirsum_{\substack{x\in L_S\\S\in\Ess}}
\Dirsum_{x\in L_S, S\in\Ess}
X_{(S,x)}
\end{displaymath}
  As a sum of metric spaces, $X$~is metrisable
  \citep[Thm.~4.2.1]{Engelking89}. Define a mapping $f\maps X\to
  Y$ by $f(y,S,x)\defeq y$. We wish to show that $f$~is a quotient
  mapping. Since every constant sequence converges, $f$~is onto.

  To show that $f$~is continuous, it is enough
  \citep[Prop.~2.2.6]{Engelking89} to show that each of the
  restrictions ${f\restriction X_{(S,x)}}$ is continuous. Let
  $S\in\Ess$, let $x\in L_S$ and consider the restricted mapping
  $g\defeq{f\restriction X_{(S,x)}}$. Let $U$~be an open subset
  of~$Y$, and let $(y,S,x)$~be a point of~$g^{-1}(U)$. Suppose that
  $\seq({(y_i,S,x)}\st i\in\nats)$ is a sequence in
  $X_{(S,x)}\setminus\set{(y,S,x)}$ that converges to~$(y,S,x)$. Then
  $y_i\neq y$ for all $i\in\nats$ and $\seq(y_i)$ converges to~$y$ in
  $S(\nats)\union\set{x}$. Now $S$~converges in~$Y$ to~$x$, hence so
  does its subsequence~$\seq(y_i)$. Since $U$~is open in~$Y$ it
  follows that $\seq(y_i)$ is eventually in~$U$. Hence
  $\seq({(y_i,S,x)})$~is eventually in $g^{-1}(U)$. This shows that
  $g^{-1}(U)$ is sequentially open, hence open, in~$X_{(S,x)}$;
  thus~$g$ is continuous. Therefore $f$~is continuous.

  So far, everything that we have done could be done for a general
  space~$Y$. Now finally we must show that $f$~is a quotient mapping:
  here we need to use the fact that $Y$~is sequential. Let $V$~be a
  subset of~$Y$ such that $f^{-1}(V)$ is open in~$X$. Let $x\in V$
  and suppose that $S=\seq(x_i\st i\in\nats)$ is a sequence in~$Y$
  that converges to~$x$. Form the subspace $X_{(S,x)}$ of~$X$, and
  note that $f^{-1}(V)\intersect X_{(S,x)}$ is a neighbourhood
  of~$(x,S,x)$ in this subspace. By definition of the subspace, the
  sequence $\seq({(x_i,S,x)})$ converges in~$X_{(S,x)}$ to~$(x,S,x)$,
  so therefore this sequence is eventually in $f^{-1}(V)\intersect
  X_{(S,x)}$. This implies that $\seq(x_i)$ is eventually in~$V$.
  Therefore $V$~is sequentially open in~$Y$, so $V$~is open because
  $Y$~is a sequential space. This shows that $f$~is a quotient
  mapping.
\end{proof}\index{sequential space|)}

%%%%%%%%%%%%%%%%%%%%%%%%%%%%%%%%%%%%%%%%%%%%%%%%%%%%%%%%%%%%%%%%%%%%%
%%%%%%%%%%%%%%%%%%%%%%%%%%%%%%%%%%%%%%%%%%%%%%%%%%%%%%%%%%%%%%%%%%%%%
%%%%                                                             %%%%
%%%%                      CHAPTER FOUR                           %%%%
%%%%                                                             %%%%
%%%%%%%%%%%%%%%%%%%%%%%%%%%%%%%%%%%%%%%%%%%%%%%%%%%%%%%%%%%%%%%%%%%%%
%%%%%%%%%%%%%%%%%%%%%%%%%%%%%%%%%%%%%%%%%%%%%%%%%%%%%%%%%%%%%%%%%%%%%

\section{Fréchet Spaces}\label{chap:frechet}\index{Fréchet space|(}

\subsection*{Unique Limits}

In this section we discuss the class of Fréchet spaces.
We have seen (Lemma~\ref{lem:firstcount-frechet}) that every first 
countable space is Fréchet. An example of a Fréchet space that is not 
first countable is the set~$\real$ with the topology
$\set{\real\setminus F\st\text{finite $F\sub\real$}}$. Shortly we are 
going to define a non-Hausdorff Fréchet space with unique sequential 
limits. This will show that the result of 
Proposition~\ref{prop:firstcount-hausdorff} does not generalise to Fréchet 
spaces. Later on (Theorem~\ref{thm:frechet-characterisation}) we are going 
to prove a characterisation theorem for the class of Hausdorff Fréchet 
spaces: the Hausdorff Fréchet spaces are precisely the images of the 
metric spaces under pseudo-open mappings. In particular, the image of a 
metric space under a closed mapping is Fréchet.

\begin{example}[{\citet[Eg.~6.2]{Franklin67}}]
\index{Hausdorff property}\label{eg:frechet-hausdorff}
Let $X$~be the product $\nats\setprod\nats$ in the usual (discrete) 
topology, and let $p$,~$q\notin X$ with $p\neq q$. Let $Y\defeq 
X\union\set{p,q}$, and define a topology on~$Y$ as follows.
The topology will be the coarsest with a base consisting of the topology 
on~$X$ and:
\begin{enumerate}
\item\label{item:nhoodp} 
the sets $\set{p}\union((\nats\setminus F)\setprod\nats)$, where $F$~is a 
finite subset of~$\nats$;
\item\label{item:nhoodq} 
the sets $\set{q}\union\Union_{i\in\nats}(\set{i}\setprod(\nats\setminus
F_i))$, where $\seq(F_i)$ is a sequence of finite subsets of~$\nats$.
\end{enumerate}
Let $U$,~$V$ be neighbourhoods of $p$,~$q$ respectively. Then $U$~contains 
all but finitely many lines $\set{i}\setprod\nats$, so we can pick a 
natural number~$i$ such that~$\set{i}\setprod\nats\sub U$. Then 
$V$~contains $(i,j)$~for all but finitely many~$j\in\nats$, so we can pick 
a natural number~$j$ such that $(i,j)$~is in~$V$. Then $(i,j)$~is in~$U$ 
too, so $U\intersect V\neq\emptyset$. This shows that every neighbourhood 
of~$p$ meets every neighbourhood of~$q$, so $Y$~is not Hausdorff.

Next we show that no sequence in~$Y$ may converge to more than one
  point. Since $X$~is discrete the only sequences converging to
  points of~$X$ are the eventually constant ones. Thus we need only
  consider the possibility that there might exist a sequence in~$Y$
  converging both to~$p$ and to~$q$. Moreover, there exist
  neighbourhoods of~$p$ that do not include~$q$, and vice versa, so we
  need only consider the case of a sequence $\seq({(i_r,j_r)}\st
  r\in\nats)$ in~$X$. Since $\seq({(i_r,j_r)})$ converges to~$p$,
  every neighbourhood of~$p$ contains all but finitely many of the
  $(i_r,j_r)$. For each $i\in\nats$, the set
  $Y\setminus(\set{i}\setprod\nats)$ is a neighbourhood of~$p$, so
  there are only finitely many $r$ such that $i_r=i$. Let
  $F_i\defeq\set{j_r\st i_r=i}$. Since $\seq(F_i)$ is a sequence of
  finite subsets of~$\nats$, the set
  \begin{displaymath}
    N\defeq\set{\smash{q}}\union\Union_{i\in\nats}(\set{i}\setprod
    (\nats\setminus F_i))
  \end{displaymath}
  is a neighbourhood of~$q$. This is a contradiction, because
  $N$~contains none of the points of the sequence $\seq({(i_r,j_r)}\st
  r\in\nats)$ that supposedly converges to~$q$. Therefore there
  exists no such sequence after all. This shows that sequential
  limits in~$Y$ are unique.

  Finally we shall show that $Y$~is a Fréchet space. Since $X$~is
  discrete, it is certainly Fréchet. It remains to prove that,
  whenever one of $p$,~$q$ lies in the closure of a set, then it is
  the limit of a sequence in that set. First consider~$p$. Suppose
  that $A$~is a subset of~$Y$ such that $p\in\closure A$. For each
  natural number~$r$, the complement of the set~$\set{i\in\nats\st
    i<r}\setprod\nats$ is a neighbourhood of~$p$. It follows that for
  each $r\in\nats$ there exists some~$(i_r,j_r)$ in~$A$ with $i_r\gee
  r$. Then $i_r\to\infty$ as $r\to\infty$. Since any neighbourhood
  of~$p$ contains almost all of the lines $\set{i}\setprod\nats$, it
  must contain almost all of the points~$(i_r,j_r)$. Thus
  $\seq({(i_r,j_r)}\st r\in\nats)$ is a sequence in~$A$ that converges
  to~$p$.

  Now consider~$q$. Let $B$~be a subset of~$Y$ with $q\in\closure B$.
  If $q\in B$, then there is nothing to prove. Suppose $q\notin B$.
  Then for each $i\in\nats$, let $F_i\defeq\set{j\in\nats\st(i,j)\in
    B}$, and consider the sequence $\seq(F_i)$ of subsets of~$\nats$.
  If every~$F_i$ is finite, then the set
  \begin{displaymath}
    \set{\smash{q}}\union(X\setminus B)
    =\set{\smash{q}}\union\Union_{i\in\nats}(\set{i}\setprod
    (\nats\setminus F_i))
  \end{displaymath}
  is of the form given in \bref{item:nhoodq}~above, so it is a
  neighbourhood of~$q$ that fails to meet~$B$, contradicting
  $q\in\closure B$. Therefore there exists a natural number~$i$ such
  that $F_i$~is infinite. Arrange the elements of~$F_i$ into an
  infinite sequence $\seq(j_r\st r\in\nats)$. Now any neighbourhood
  of~$q$ contains almost all points of the
  line~$\set{i}\setprod\nats$, in particular it contains almost all of
  the $(i,j_r)$. Thus $\seq({(i,j_r)}\st r\in\nats)$ is a sequence
  in~$B$ that converges to~$q$. This shows that $Y$~is a Fréchet
  space.
\qed
\end{example}

It is possible to give a slightly simpler example of a
non-Hausdorff Fréchet space with unique sequential limits
\citep[]{Franklin66}. Note however that in
Example~\ref{eg:frechet-hausdorff}, the space $Y$~is compact. This is
of interest, because recalling Example~\ref{eg:sequential-productive}
we can deduce that the square of
a compact \index{product!of Fréchet space}%
\index{product!of sequential space}%
Fréchet space need not be sequential. In fact, the situation is
about as bad as it could possibly be: even the product of a
Hausdorff Fréchet space with a compact metric space need not be
Fréchet. For instance, if~$X$ is the quotient space of~$\real$
obtained by identifying the integers to a point, then it is an
easy exercise to show that $X$~is a Fréchet space but
$X\setprod\I$ is not (this example is due to
\citet[Eg.~7.4]{Franklin67}). Products of spaces of countable
tightness are a little better behaved: see
Theorem~\ref{thm:counttight-compactproductive}.

Like first countability and unlike sequentiality, the Fréchet
property is hereditary:
\begin{lem}\label{lem:frechet-hereditary}\index{subspace!of Fréchet space}
  Any subspace of a Fréchet space is itself a Fréchet space.\qed
\end{lem}
This result furnishes a converse to
Lemma~\ref{lem:hereditarilysequential-frechet}; so we can deduce
that the Fréchet spaces are precisely the hereditarily sequential
spaces.

\subsection*{Characterisation Theorem}

\begin{figure}[tp]
\begin{center}
\small
\begin{displaymath}
%%% DIAGRAM: needs xymatrix package installed \usepackage[all]{xy}
\renewcommand{\objectstyle}{\text}
\xymatrix{%
			&			&proper \ar[d]		\\
open onto \ar[dr]	& 			&closed onto \ar[dl]	\\
			&pseudo-open \ar[d]	& 			\\
			&quotient		&
}
\end{displaymath}
\caption{Relations between some classes of continuous mappings.}
\label{fig:mappingrelations}
\end{center}
\end{figure}
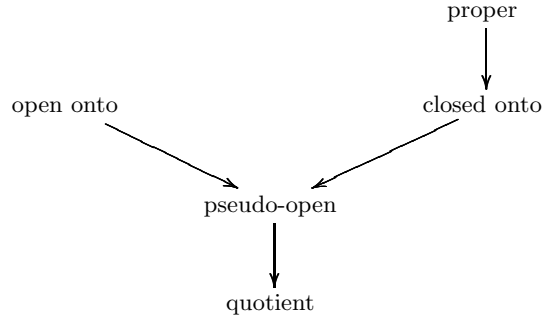

Before we state our characterisation theorem for Fréchet spaces,
we need to define a new class of mappings.
\begin{defn}[pseudo-open mapping]
  Let $X$,~$Y$ be Hausdorff spaces, and let $f\maps X\to Y$ be
  continuous and onto. We say that $f$~is \define[pseudo-open
  mapping]{pseudo-open} if whenever $U$~is an open set in~$X$
  containing some fibre $f^{-1}\set{y}$, then $f(U)$~is a
  neighbourhood of~$y$ in~$Y$.
\end{defn}
It is clear that every open onto mapping is pseudo-open. Also, any
closed onto mapping is pseudo-open: for if $f\maps X\to Y$ is
closed, then let $y\in Y$. Suppose that $U$~is open in~$X$ with
$f^{-1}\set{y}\sub U$. Then $X\setminus U$ is closed in~$X$, so
$f(X\setminus U)$ is closed in~$Y$. Since $f^{-1}\set{y}$~is
contained in~$U$, it follows that $y$~is not an element
of~$f(X\setminus U)$. Hence $y$~belongs to the open set
$Y\setminus f(X\setminus U)$, which is a subset of~$f(U)$ (because
$f$~is onto). Thus $f(U)$~is a neighbourhood of~$y$ and $f$~is
pseudo-open. Figure~\ref{fig:mappingrelations} summarises the
situation.

%\medskip
%\pagebreak[3]

Now we can present the characterisation theorem. This theorem is
due to \citet[Thm.~2]{Arhangelskij63}, although in our proof we
follow \citet{Franklin65}.
\begin{thm}[Arhangel{\textprime}ski{\u\i}]
  \label{thm:frechet-characterisation}\index{pseudo-open mapping}
  If $f$~is a pseudo-open mapping from a metric space~$X$ onto a
  Hausdorff space~$Y$, then $Y$~is Fréchet. Conversely, if\/ $Y$~is a
  Hausdorff Fréchet space, then for some metric space~$X$ there exists
  a pseudo-open mapping~$f$ from~$X$ onto~$Y$.
\end{thm}

To prove Theorem~\ref{thm:frechet-characterisation} we shall call upon
Theorem~\ref{thm:sequential-characterisation}. First, we note that
every pseudo-open mapping is quotient:
\begin{lem}\label{lem:pseudoopen-quotient}
  \index{quotient mapping}\index{pseudo-open mapping}
  Let $X$,~$Y$ be Hausdorff spaces, and let $f\maps X\to Y$ be
  pseudo-open. Then $f$~is a quotient mapping.\qed
\end{lem}

The next result will allow us to give an easy proof of
Theorem~\ref{thm:frechet-characterisation}.
\begin{prop}[{\citet[Thm.~4]{Arhangelskij63}}]
  \label{prop:quotient-pseudoopen}\index{quotient mapping}
  \index{pseudo-open mapping}
  Let $X$~be a Hausdorff Fréchet space, let $Y$~be a Hausdorff space
  and let $f\maps X\to Y$ be a quotient map. Then $Y$~is a Fréchet
  space if and only if $f$~is pseudo-open.
\end{prop}
\begin{proof}(Franklin)
  Suppose first that $Y$~is a Fréchet space. Let $U$~be an open
  subset of~$X$ containing some fibre $f^{-1}\set{y}$. Suppose for a
  contradiction that $f(U)$~is not a neighbourhood of~$y$. Then
  $y\in\closure{(Y\setminus f(U))}$, so as $Y$~is Fréchet there exists
  a sequence $\seq(y_i)$ in~$Y\setminus f(U)$ that converges to~$y$.
  Let $A\defeq\set{y_i\st i\in\nats}$, and let $B\defeq f^{-1}(A)$.
  Note that $\closure A=A\union\set{y}$ because $Y$~is Hausdorff and
  $\seq(y_i)$ converges to~$y$. Therefore, since $f$~is continuous,
  \begin{equation}\label{eq:3}
    \closure B=\closure{(f^{-1}(A))}\sub f^{-1}(\closure A)
    =f^{-1}(A\union\set{\smash{y}})=B\union f^{-1}\set{\smash{y}}
  \end{equation}
  However, we have $A\sub Y\setminus f(U)$, whence $U\intersect
  B=\emptyset$. Since $U$~is open, it follows that
  $U\intersect\closure B=\emptyset$. Now we recall that
  $f^{-1}\set{y}\sub U$, so therefore $f^{-1}\set{y}\intersect
  \closure B=\emptyset$. From~\bref{eq:3} we now see that $\closure
  B\sub B$, i.e., that $B$~is closed. Therefore $f^{-1}(Y\setminus
  A)=X\setminus B$ is open, and since $f$~is a quotient mapping it
  follows that $Y\setminus A$~is open, i.e., that $A$ is closed. Thus
  $A=\closure A=A\union\set{y}$, so that $y\in A\sub Y\setminus f(U)$.
  This contradicts the fact that $f^{-1}\set{y}\sub U$. Therefore
  $f(U)$~is a neighbourhood of~$y$ after all, hence $f$~is
  pseudo-open.

  Now suppose conversely that $f$~is pseudo-open: we wish to show that
  $Y$~is Fréchet. Let $A$~be a subset of~$Y$ and suppose that
  $y\in\closure A$. Let $B\defeq f^{-1}(A)$ and suppose for a
  contradiction that $f^{-1}\set{y}\intersect\closure B=\emptyset$.
  Then there exists an open set~$U$ containing~$f^{-1}\set{y}$ such
  that $U\intersect B=\emptyset$, whence $f(U)\intersect A=\emptyset$.
  Since $f$~is pseudo-open, $f(U)$~is a neighbourhood of~$y$, and this
  contradicts the fact that $y\in\closure A$. Therefore
  $f^{-1}\set{y}\intersect\closure B\neq\emptyset$, so there exists
  some $x\in f^{-1}\set{y}$ such that $x\in\closure B$. Since $X$~is
  a Fréchet space, there exists a sequence $\seq(x_i)$ in~$B$ that
  converges to~$x$. Then $\seq({f(x_i)})$ is a sequence in~$A$ that
  converges to $f(x)= y$. This shows that $Y$~is a Fréchet space.
\end{proof}
\begin{cor}\index{proper mapping!of Fréchet space}
  \label{cor:frechet-properimage}
  The class of Fréchet spaces is closed with respect to taking
  Hausdorff images under proper mappings.\qed
\end{cor}

At last we arrive at the proof of the main theorem.
\begin{proof}[Proof of\/ Theorem~\ref{thm:frechet-characterisation}]
  Suppose first that $X$~is a metric space, that $Y$~is a Hausdorff
  space and that $f\maps X\to Y$ is pseudo-open and onto. Since
  $X$~is a metric space, it is Hausdorff and Fréchet. Therefore by
  Proposition~\ref{prop:quotient-pseudoopen}, \ $Y$~is a Fréchet space.

  Conversely, suppose that $Y$~is a Hausdorff Fréchet space. In
  particular $Y$~is a sequential space, so by
  Theorem~\ref{thm:sequential-characterisation} there exists a metric
  space~$X$ and a quotient mapping $f\maps X\to Y$. As a metric
  space, $X$~is Hausdorff and Fréchet, therefore $f$~is pseudo-open by
  Proposition~\ref{prop:quotient-pseudoopen}.
\end{proof}
\index{Fréchet space|)}

%%%%%%%%%%%%%%%%%%%%%%%%%%%%%%%%%%%%%%%%%%%%%%%%%%%%%%%%%%%%%%%%%%%%%
%%%%%%%%%%%%%%%%%%%%%%%%%%%%%%%%%%%%%%%%%%%%%%%%%%%%%%%%%%%%%%%%%%%%%
%%%%                                                             %%%%
%%%%                      CHAPTER FIVE                           %%%%
%%%%                                                             %%%%
%%%%%%%%%%%%%%%%%%%%%%%%%%%%%%%%%%%%%%%%%%%%%%%%%%%%%%%%%%%%%%%%%%%%%
%%%%%%%%%%%%%%%%%%%%%%%%%%%%%%%%%%%%%%%%%%%%%%%%%%%%%%%%%%%%%%%%%%%%%

\section{Countable Tightness}\index{countable tightness|(}
\label{chap:countable-tightness}

\subsection*{Preamble}

Countable tightness is of interest to us principally because of the
Moore-Mrówka problem. The solution of the problem has been described as
`the main advance in the theory of compact spaces during [the 1980s],' by
\citet[p.~574]{Shakhmatov92}. This indicates the significance of tightness
and sequentiality in analytic topology. One of the trends in topology
since the 1960s has been the growing importance of \index{cardinal
invariants}% cardinal invariants of spaces, of which tightness is a
typical example. These invariants are usually defined in general terms,
but most of them are generalisations of preexisting concepts involving
countability in some way. A summary of relations between some of the most
important cardinal invariants appears in
Figure~\ref{fig:cardinal-relations}. For more details (and proofs) see
\citet{Hodel84} and \citet{Juhasz80,Juhasz84}.

\begin{figure}[h]
\begin{minipage}{\linewidth}
\centering
%%% DIAGRAM: needs xymatrix package installed \usepackage[all]{xy}
\begin{displaymath}
\renewcommand{\labelstyle}{\textstyle}
\xymatrix{%
	&\abs{X} \ar[d]			&w(X) \ar[dl] \ar[dd]			&			\\
	&nw(X) \ar[d] \ar[drr]		&					&			\\
	&hL(X) \ar[dl] \ar[d]^{T_3}	&\chi(X) \ar[ddl]^(.6){T_1} \ar[dd]	&hd(X) \ar[d] \ar[ddl]	\\
L(X)	&\Psi(X) \ar[d]_{T_1}		&					&d(X) \ar[d]		\\
	&\psi(X)			&t(X)					&c(X)			\\
}
\end{displaymath}
\index{cardinal invariants}
\index{countable tightness}
\index{countable pseudo-character}
\index{index of perfection}
\index{hereditarily Lindelöf space}
\index{hereditarily separable space}

%\bigskip

\subsection*{Key}
\begin{tabular}{lll}
Invariant	&Name					&If $\text{invariant}=\aleph_0$, space is called	\\
\hline
$\abs{X}$	&cardinality of $X$			&countable						\\
$w(X)$		&weight of $X$				&second countable					\\
$nw(X)$		&network weight of $X$			&							\\
$\chi(X)$	&character of $X$			&first countable					\\
$hL(X)$		&hereditary Lindelöf number of $X$	&hereditarily Lindelöf					\\
$hd(X)$		&hereditary density of $X$		&hereditarily separable					\\
$\Psi(X)$	&index of perfection\footnote{This 
		terminology is non-standard. See 
		Section~\ref{chap:perfect}.} of $X$	&perfect						\\
$L(X)$		&Lindelöf number of $X$			&Lindelöf						\\
$d(X)$		&density of $X$				&separable						\\
$t(X)$		&tightness of $X$			&							\\
$\psi(X)$	&pseudo-character of $X$		&							\\
$c(X)$		&Souslin number of $X$			&CCC (countable chain condition)			\\
\end{tabular}
\end{minipage}
\caption{Relations between cardinal invariants.}
\label{fig:cardinal-relations}
\end{figure}

\subsection*{How about Sequences?}

We noted in Section~\ref{chap:intro} that every sequential space has
countable tightness. The converse is false, however.
\begin{example}[{\citet[Thm.~A]{Franklin69}}]
  \label{eg:counttight-sequential}\index{sequential space}
  We exhibit a non-sequential space with countable tightness. Let
  $\beta\nats$ be the Stone-{\v C}ech compactification of the natural
  numbers. Let $X\defeq\beta\nats$ and consider the topology on~$X$
  generated by the topology on~$\beta\nats$ together with the sets
  $\set{x}\union\nats$ for all $x$ in $X\setminus\nats$. We have
  already noted that the only convergent sequences in~$\beta\nats$ are
  those that are eventually constant (Example~\ref{eg:betanats}).
  Clearly, in a finer topology the number of convergent sequences
  cannot increase, so the only convergent sequences in~$X$ are those
  that are eventually constant. Therefore every subset of~$X$ is
  sequentially open, but $X\setminus\nats$ is not open, so $X$~is not a
  sequential space. However, $X$~has countable tightness. To see
  this, let $A$~be any subset of~$X$. Note that any subset of~$X$
  containing~$\nats$ is open, hence any subset that is disjoint
  from~$\nats$ is closed. Therefore
  \begin{equation}\label{eq:4}
    \closure{A}
    =\closure{(A\setminus\nats)}\union\closure{(A\intersect\nats)}
    =(A\setminus\nats)\union\closure{(A\intersect\nats)}
    \sub A\union\closure{(A\intersect\nats)}\text.
  \end{equation}
  Now, if $A$~contains the closure of all its countable subsets then
  in particular it contains the closure of $A\intersect\nats$. Then
  \bref{eq:4}~shows that $\closure A=A$, so that $A$~is closed. This
  shows that $X$~has countable tightness.
\qed
\end{example}

Using the last example, we can show that there is no
generalisation of Lemma~\ref{lem:sequential-conts} to spaces of
countable tightness.
\begin{example}\index{continuous mapping}
  \label{eg:counttight-conts}
  Let $X$~be the space of Example~\ref{eg:counttight-sequential}, and let
  $Y$~be the same set ($\beta\nats$) in the discrete topology. Then
  the topology on~$Y$ is strictly finer than the topology on~$X$, so
  the identity mapping $\iota\maps X\to Y$ is \emph{not} continuous.
  However, the only convergent sequences in~$X$ are those that are
  eventually constant, hence $\iota$~preserves limits of sequences.
\qed
\end{example}

With this example and previous results, all the assertions of
Figure~\ref{fig:relations} are proved.

\subsection*{Countable Tightness is not Productive}

How well does the class of spaces of countable tightness behave
under the usual operations?  First of all we note the elementary
result that countable tightness is hereditary.
\begin{lem}\index{subspace!of countably tight space}
  \label{lem:counttight-hereditary}
  Let $X$~be a topological space with countable tightness, and let
  $Y$~be a subspace of\/~$X$. Then $Y$~has countable tightness.\qed
\end{lem}
To show that countable tightness is an invariant of proper
mappings, we use the following:
\begin{lem}[{\citet[Prop.~3]{Arhangelskij68}}]
  \index{quotient mapping!of countably tight space}
  \label{lem:counttight-quotient}
  Let $X$~be a topological space of countable tightness, and suppose
  that $f\maps X\to Y$ is a quotient mapping. Then $Y$~has countable
  tightness.
\end{lem}
\begin{proof}
  Let $A$~be a subset of~$Y$ that contains the closure of all its
  countable subsets. Let $C$~be a countable subset of $f^{-1}(A)$.
  Then $f(C)$~is a countable subset of~$A$, so $\closure{f(C)}\sub A$.
  Using the fact that $f$~is continuous, it follows that
  \begin{equation}\label{eq:5}
    \closure{C}\sub\closure{\bigl(f^{-1}f(C)\bigr)}\sub
      f^{-1}\bigl(\closure{f(C)}\bigr)\sub f^{-1}(A)\text.
  \end{equation}
  Therefore $f^{-1}(A)$ is a subset of~$X$ that contains the closure
  of each of its countable subsets. Since $X$~has countable
  tightness, we can deduce that $f^{-1}(A)$~is closed in~$X$. Since
  $f$~is a quotient mapping, it follows that $A$~is closed in~$Y$.
  Hence $Y$~has countable tightness.
\end{proof}
\begin{cor}
  \index{proper mapping!of countably tight space}
  \label{cor:counttight-properimage}
  The class of spaces with countable tightness is closed with respect
  to taking images under proper mappings.\qed
\end{cor}

The product of two spaces of countable
\index{product!of countably tight space}%
tightness need not in general have countable tightness. An example
was given by \citet{Arhangelskij72}: we omit the details. Let
$X$~be the quotient space obtained by identifying the limit points
in the disjoint sum of countably many convergent sequences, and
let $Y$~be the corresponding space constructed from uncountably
many convergent sequences. Then it is easily verified that~$X$
and~$Y$ have countable tightness. However, their product
$X\setprod Y$ has uncountable tightness (consider the point whose
coordinates are the limit points in~$X$,~$Y$ respectively).

If one of the factor spaces is compact then the product
\emph{does} have countable tightness, according to the following
theorem. The proof we give is due to \citet[5.9]{Juhasz80}.
\begin{thm}[{\citet[Thm.~4]{Malyhin72}}]
  \label{thm:counttight-compactproductive}
  \index{product!of countably tight space}\index{compact space}
  Let $X$,~$Y$ be Hausdorff spaces with countable tightness, and
  suppose that $Y$~is compact. Then $X\setprod Y$~has countable
  tightness.
\end{thm}
\begin{proof}
  Let $A$~be a subset of~$X\setprod Y$ that contains the closure of
  all its countable subsets. We wish to show that $A$~is closed.
  Suppose for a contradiction that it is not, so that there exists a
  point $(x,y)$ in $\closure A\setminus A$. It is not obvious how to
  obtain a contradiction, so a few remarks are in order before we
  proceed. We know very little about the product $X\setprod Y$, so
  most of the work will have to be done in the factor spaces. The
  plan of attack is as follows: we are going to define two subsets
  $B$~and~$W$ of~$A$ in such a way that they are obviously disjoint,
%  (see Figure~\ref{fig:counttight-prod}), 
  and then we shall show that
  they have a point in common. More specifically, we define~$B$ and
  project it down to~$Y$, then we use the separation properties of~$Y$
  and the assumption ${(x,y)\notin A}$ to define~$W$ disjoint
  from~$B$. Then we project $W$~down to~$X$ and show that the
  assumption ${(x,y)\in\closure A}$ implies that $W$~and~$B$
  \emph{cannot} be disjoint after all.

  Since $X$~is a \mbox{$T_1$-space}, $\set{x}$~is closed in~$X$, hence
  $\set{x}\setprod Y$ is closed in $X\setprod Y$. The subspace
  $\set{x}\setprod Y$~is homeomorphic to~$Y$, so it has countable
  tightness. Let $B\defeq A\intersect(\set{x}\setprod Y)$. Then
  $B$~contains the closure in $\set{x}\setprod Y$ of all its countable
  subsets. Since $\set{x}\setprod Y$ has countable tightness, it
  follows that $B$~is closed in~$\set{x}\setprod Y$, hence
  in~$X\setprod Y$. Now, $\set{x}\setprod Y$ is homeomorphic to~$Y$,
  and $B$~is closed in~$\set{x}\setprod Y$, so therefore $\pi_Y(B)$ is
  closed in~$Y$. Since $(x,y)\notin A$ we have $y\notin\pi_Y(B)$.
  Now $Y$~is a regular space, so there exist disjoint open subsets
  $U$,~$V$ of~$Y$ such that $y\in U$ and $\pi_Y(B)\sub V$. So if we
  let $E\defeq\closure[Y] U$, then $E$~is a closed neighbourhood
  of~$y$ in~$Y$ and $\pi_Y(B)\intersect E=\emptyset$. Hence
  $X\setprod E$ is a closed neighbourhood of~$(x,y)$ in~$X\setprod Y$.
%    \begin{figure}[ht]
%    \begin{center}
%      \includegraphics{countabletightness.eps}
%    \end{center}
%    \caption{}
%    \label{fig:counttight-prod}
%  \end{figure}

  Finally now, let $W\defeq A\intersect(X\setprod E)$. Let $N$~be a
  neighbourhood of~$(x,y)$. Then ${(X\setprod E)}\intersect N$~is
  another neighbourhood, and since $(x,y)$~belongs to~$\closure A$,
  the set $W\intersect N=A\intersect{(X\setprod E)}\intersect N$ is
  non-empty. This shows that $(x,y)$~is an element of the closure
  of~$W$. Now $W$~contains the closure of all its countable subsets,
  because it is the intersection of two sets with this property.
  Since $Y$~is compact, the projection~$\pi_X$ is closed. Hence
  $\pi_X(W)$~is a subset of~$X$ that contains the closure of all its
  countable subsets. Since $X$~has countable tightness, it follows
  that $\pi_X(W)$~is closed in~$X$. Using the continuity of~$\pi_X$,
  we have
  \begin{equation}\label{eq:2}
    \pi_X(W)=\closure[X]{\pi_X(W)}\contains\pi_X(\closure W)\ni x\text.
  \end{equation}
  Hence there exists some point~$(x,z)$ in $W=A\intersect{(X\setprod
    E)}$. Then $z$~belongs to~$E$, and $(x,z)$~belongs
  to~$A\intersect(\set{x}\setprod Y)=B$. Thus we arrive at the
  contradiction that $z$~is an element of the set $\pi_Y(B)\intersect
  E$, a set which is empty by our choice of~$E$.
\end{proof}

\subsection*{The Moore-Mrówka Problem}

It is worth remarking that in Example~\ref{eg:counttight-sequential}
we had to adapt the topology of the non-sequential
space~$\beta\nats$ in such a way that it was no longer compact,
before we obtained countable tightness. The question
\begin{enumerate}\index{sequential space}\index{compact space}
\item[(MM)] Is every compact Hausdorff space of countable
tightness a
  sequential space?
\end{enumerate}
is known as the \define{Moore-Mrówka problem}. The problem was
posed in \citet{Moore64}. It is now known that the question cannot
be answered in ZFC\@. The set-theory required to prove this means
that it would take too long to go into the details here, so we
must content ourselves with stating the following two-part
theorem:

%\pagebreak[4]
\begin{thm}\label{thm:moore-mrowka}
%\par\indent
\mbox{}
\begin{enumerate}
\item 
Assume~$\diamondsuit$, and let $X$~be Ostaszewski's space 
\citep[]{Ostaszewski76}. Let $\alpha X$~be the Alexandroff one-point 
compactification of\/~$X$. Then $\alpha X$~has countable tightness but is 
not a sequential space.
\item 
Assume~PFA\@. Then every compact Hausdorff space of countable
tightness is a sequential space---see \citet{Balogh89}.
\end{enumerate}
\end{thm}

Many other results connected to the Moore-Mrówka problem can be found in 
the survey by \citet[Sec.~2]{Shakhmatov92}.\index{countable tightness|)}

%%%%%%%%%%%%%%%%%%%%%%%%%%%%%%%%%%%%%%%%%%%%%%%%%%%%%%%%%%%%%%%%%%%%%
%%%%%%%%%%%%%%%%%%%%%%%%%%%%%%%%%%%%%%%%%%%%%%%%%%%%%%%%%%%%%%%%%%%%%
%%%%                                                             %%%%
%%%%                      CHAPTER SIX                            %%%%
%%%%                                                             %%%%
%%%%%%%%%%%%%%%%%%%%%%%%%%%%%%%%%%%%%%%%%%%%%%%%%%%%%%%%%%%%%%%%%%%%%
%%%%%%%%%%%%%%%%%%%%%%%%%%%%%%%%%%%%%%%%%%%%%%%%%%%%%%%%%%%%%%%%%%%%%

\section{Perfect Spaces}\label{chap:perfect}\index{perfect space|(}

\subsection*{Outline}

We are now going to consider perfect spaces: we are particularly
interested in compactness in perfect spaces. We start off by considering
the relationship between the class of perfect spaces and other important
classes of spaces. In general, the perfect spaces are not directly related
to those that we have considered up until now
(Examples~\ref{eg:firstcount-perfect} and~\ref{eg:perfect-counttight}).
However, when suitable compactness properties and separation axioms are
assumed, then some useful results are achieved. In particular, we are
going to prove that every countably compact regular perfect space has
countable tightness (Theorem~\ref{thm:perfectcompact-counttight}).

The reader is asked to look again at Figure~\ref{fig:relations}. Of
the many useful properties of sequences in metric spaces, we have
generalised all except the following: `compactness and sequential
compactness are equivalent'. The statement is true for first
countable paracompact spaces. However, paracompactness is a very
strong condition. Looking to strengthen the result, it is
noteworthy that every perfectly normal space is \emph{countably}
paracompact \citep[Cor.~to~V.5]{Nagata85}. Since every metric
space is both first countable and perfectly normal, an attractive
conjecture is that compactness and sequential compactness coincide
for these spaces---in fact this is independent of ZFC
(Corollary~\ref{cor:perfectfirstcount-compact}).

\subsection*{When is a Space Perfect?}

The following result appears in \citet[Chap.~3, Prob.~206]{Arhangelskij84a}.

\begin{lem}\index{hereditarily Lindelöf space}\label{lem:heredlindelof-perfect}
Suppose that $X$~is hereditarily Lindelöf and regular. Then $X$~is a 
perfect space.
\end{lem}

\begin{proof}
  Let $U$~be any open subset of~$X$. Since $X$~is regular, for every
  $x\in U$ there exists an open set~$V$ such that $x\in V\sub\closure
  V\sub U$. Then the family consisting of all the sets $V$ is an open
  cover for~$U$. Since $X$~is hereditarily Lindelöf, $U$~is Lindelöf.
  So there exists a countable subcover $\set{V_i\st i\in\nats}$. But
  by construction $\closure{V_i}\sub U$ for every natural number~$i$,
  and since $\set{V_i}$ is a cover, it follows that
  $U=\Union_{i\in\nats}\closure{V_i}$. This shows that $U$~is an
  $F_\sigma$. Hence $X$~is a perfect space.
\end{proof}

\begin{cor}\index{second axiom of countability}
  \index{second countable|see{second axiom of countability}}
  \label{cor:secondcount-perfect}
  A second countable regular space is perfect.\qed
\end{cor}

Of course, by Urysohn's metrisation theorem a second countable
regular space is pseudo-metrisable, and pseudo-metric spaces are
always perfect. However, our argument generalises to arbitrary
cardinalities. Thus for regular spaces, $\Psi(X)\lee hL(X)$, where
$hL(X)$~is the hereditary Lindelöf number of the space, and
\label{perfectionindex}%
\define[psix@$\Psi(X)$|see{index of perfection}]{$\Psi(X)$}~%
is the \define{index of perfection}, %
that is, the least infinite cardinal~$\cardm$ such that every
closed set can be written as an intersection of~${}\lee\cardm$
open sets.

It is well known that every regular Lindelöf space is normal
(\citet[Lem.~4.1]{Kelley55} calls this
\define[Tychonoff's lemma]%
{\mbox{Tychonoff's} lemma}). Thus
Lemma~\ref{lem:heredlindelof-perfect} implies that a regular
hereditarily Lindelöf space is perfectly normal.

It should be noted in this context that there exist first
countable compact Hausdorff spaces that are not perfect:

\begin{example}[{Alexandroff double circle \citep[]{Alexandroff29}}]
\label{eg:firstcount-perfect}
Consider the two concentric circles in $\real^2$ defined by
$C_1\defeq\set{(r\angle\theta)\in\real^2\st r^2=1}$, and
$C_2\defeq\set{(r\angle\theta)\in\real^2\st r^2=2}$. For each
natural number~$i$ and each point $z=(1\angle\theta)$ in $C_1$, define two 
points $z_{i+}$ and $z_{i-}$ by $z_{i\pm}\defeq(1\angle\theta\pm 1/i)$. 
Define a map $p\maps C_1\to C_2\maps z\mapsto\sqrt{2}z$, so that $p$~is 
the projection of~$C_1$ onto $C_2$ through~$0$. The \define{Alexandroff 
double circle} is the set $X\defeq C_1\union C_2$, with the topology 
on~$X$ generated by all the sets $\set{z}$ for~$z$ in~$C_2$ and the sets
\begin{displaymath}
U_i(z)
\defeq V_i(z) \union p(V_i(z)\setminus\set{z})\text,
\qquad
\text{for $z\in C_1$ and $i\in\nats$,}
\end{displaymath}
where $V_i(z)$ is the open arc of $C_1$~between~$z_{i-}$ and~$z_{i+}$, 
i.e.,
\begin{displaymath}
V_i(r\angle\theta) 
\defeq \set{(r\angle\theta+\delta)\st -\reciprocal i<\delta<\reciprocal i}.
\end{displaymath}
It is clear that~$X$, so topologised, is a first countable Hausdorff 
space.

We shall show that $X$~is compact. Let $\Eee$ be an open cover for~$X$,
consisting of open sets from the subbasis described above. Now,
$f\maps[0,2\pi]\to C_1\maps\theta\mapsto(1\angle\theta)$ is a continuous
mapping of~$[0,2\pi]$ onto~$C_1$ in the subspace topology. Hence $C_1$~is
a compact subspace of~$X$. So there exists a finite subset
$\set{W_1,\dots,W_n}$ of~$\Eee$ that covers~$C_1$. Then as the elements
of~$\Eee$ are subbasic, each~$W_k$ is of the form $U_{i(k)}(z_k)$ for some
$i_k\in\nats$ and $z_k\in C_1$. Clearly then, $\set{W_1,\dots,W_n}$ covers
all of~$C_2$ except for possibly the points $p(z_1)$, $p(z_2)$,
$\dots$,~$p(z_n)$. For each $k=1$, $2$, \dots,~$n$, let $Z_k$~be an
element of~$\Eee$ such that $p(z_k)$~belongs to~$Z_k$. Then
$\set{W_k,Z_k\st k=1,\dots,n}$ is a finite subcover of~$\Eee$. Thus every
cover of~$\Eee$ consisting of subbasic open sets has a finite subcover. By
the Alexander theorem \citep[Thm.~5.6]{Kelley55} this is enough to show
that $X$~is a compact space. Of course, the Alexander theorem is not
required to prove this, but it does simplify the proof.

Finally we show that $X$~is not a perfect space. The
set~$C_1=X\setminus\set{\set{z}\st z\in C_2}$ is closed, but it is not a
$G_\delta$~set. For if $W$~is an open neighbourhood of~$C_1$, then it must
contain all but finitely many points of~$C_2$. Hence any countable
intersection of neighbourhoods of~$C_1$ meets~$C_2$, and it follows that
$C_1$~is not a~$G_\delta$. This shows that $X$~is a compact Hausdorff
first countable space that is not perfect.
\qed
\end{example}

What closure properties are possessed by the class of perfect spaces? The 
inverse image of a perfect space under a proper mapping need not be 
perfect, as we see if we map the Alexandroff double circle to a single 
point. On the other hand, the following results can be proved in the 
obvious fashion:

\begin{lem}\label{lem:perfect-hereditary}\par\indent
\begin{enumerate}
\item\index{proper mapping!of perfect space} 
The image of a perfect space under a closed onto mapping is again a 
perfect space. In particular, the class of perfect spaces is closed with 
respect to taking images under proper mappings.
\item\index{subspace!of perfect space} 
Every subspace of a perfect space is perfect.
\qed
\end{enumerate}
\end{lem}

We shall complete the task of showing that perfect spaces do not fit into 
the hierarchy of Figure~\ref{fig:relations}. Assuming CH, 
\citet[Thm.~1]{Hajnal74} gave an example of an \mbox{$L$-space} in which 
every countable subset is closed.\footnote{An \emph{$L$-space} is a 
regular space that is hereditarily Lindelöf but not hereditarily 
separable.} It follows that their space is perfectly normal with 
uncountable tightness. I posed the question of whether there exists such a 
space in ZFC; it cannot be compact, as we shall see  
in~Theorem~\ref{thm:perfectcompact-counttight}. The \citet{Arens50} 
example of a perfectly normal space that is not sequential cannot be used 
here, because it is countable, hence automatically has countable 
tightness. Jo Lo has supplied me with the following example, which is a 
modification of the `Radial Interval' topology described in 
\citet[141]{Steen96}. With hindsight it appears that Lo's example and 
Arens's are actually quite similar.

\begin{example}\index{countable tightness}\label{eg:perfect-counttight}
Let $\I$ denote the unit interval in the usual topology, and let $L$~be
the disjoint topological sum $\Dirsum_{r\in\real}\I\setprod\set{r}$. Let
$M$~be the quotient space of~$L$ in which all the points~$(0,r)$ are
identified (the resulting point is henceforth denoted~$0$). Let $q\maps
L\to M$ be the canonical mapping. Let $X$~be the set~$M$ with the topology
whose closed sets are the closed sets in~$M$ together with all countable
unions of sets of the form $q((0,1]\setprod\set{r})$. Then $X$~does not
have countable tightness. For, the set $X\setminus\set{0}$ is not closed,
but $0\notin\closure C$ for any countable $C\sub X\setminus\set{0}$. To
see this, note that for at most countably many values of~$r$ can the set
$C\intersect q((0,1]\setprod\set{r})$ be non-empty. If these values of~$r$
are enumerated $\set{r_i\st i\in\nats}$, then we have
\begin{equation}\label{eq:6}
\closure C\sub \Union_{i\in\nats}q((0,1]\setprod\set{r_i})\text,
\end{equation}
so that $0\notin\closure C$. Thus $X$~does not have countable
tightness.

Now we show that $X$~is a perfectly normal space. First, note that $L$~is
a topological sum of metric spaces, so it is metric. The mapping $q$~is
closed, because $q^{-1}\set{0}=\set{(0,r)\st r\in\real}$ is closed in~$L$.
It follows that $M=q(L)$ is perfectly normal; combine
Lemma~\ref{lem:perfect-hereditary}(1) with the elementary fact that closed
mappings preserve normality. Now, $X$~has mostly the same closed sets
as~$M$, and the extra ones of the form $q((0,1]\setprod\set{r})$ are both
open and closed. It is easy to see that perfect normality is preserved
under such an extension of the topology. So $X$~is perfectly normal with
uncountable tightness, as claimed.
\qed
\end{example}

\subsection*{On Compactness}\index{compact space|(}

Now we begin our consideration of compactness in perfect spaces. First we
are going to see how its presence affects matters, and later on we shall
consider when its presence may be guaranteed by weaker covering
properties.

We are going to start by proving the next theorem, which is noted in 
\citet[Prop.~3]{Weiss78}. It is used in a crucial way by \citet{Weiss78} 
in proving Theorem~\ref{thm:perfectcountcompact-compact}(2) below.

\begin{thm}\index{countable tightness}\index{countably compact space}
\label{thm:perfectcompact-counttight}
Every perfectly regular countably compact space has countable tightness.
\end{thm}

We shall prove this using two subsidiary results. The first is a 
proposition of \citet[Prop.~2.1]{Ostaszewski76}, which states that perfect 
countably compact spaces have \define{countable spread}.

\begin{prop}\label{prop:perfectcompact-countspread}
Let $X$~be a perfect countably compact space. Whenever $Y$~is a subset 
of\/~$X$ that is discrete in the subspace topology, then $Y$~is countable.
\end{prop}

\begin{proof}\index{omega@$\omega_1$}
  Suppose for a contradiction that $Y$~is uncountable. Without loss
  of generality we may assume that $\abs{Y}=\aleph_1$, so we can index
  $Y$ by the countable ordinal numbers:
  $Y=\set{x_\alpha\st\alpha<\omega_1}$. Since $Y$~is discrete, for
  every~$\alpha<\omega_1$ we can choose an open subset $U_\alpha$
  of~$X$ such that $U_\alpha\intersect Y=\set{x_\alpha}$. Now let
  $U\defeq\Union\set{U_\alpha\st\alpha<\omega_1}$. Then $U$~is an
  open set, and since $X$~is perfect it follows that~$U$ is an
  $F_\sigma$. Let $\seq(F_i)$~be an increasing sequence of closed
  sets such that $U=\Union_{i\in\nats}F_i$. Then
  \begin{equation}\label{eq:7}
    Y=U\intersect Y
    =\Union_{i\in\nats}(F_i\intersect Y)\text.
  \end{equation}
  Since $Y$~is uncountable, there is some natural number~$i$ such that
  $F_i\intersect Y$~is uncountable. Let $\alpha<\omega_1$ be such
  that $x_\alpha$~belongs to~$F_i$. Then the set
  $\set{x_{\beta}\st\beta\lee\alpha}$ is countable, so it does not
  exhaust~$F_i\intersect Y$. This fact allows us to pick a strictly
  increasing sequence $({\alpha(r)})$ of ordinal numbers~${}<\omega_1$
  such that $x_{\alpha(r)}$~belongs to~$F_i$ for all $r\in\nats$.

  The rest of the proof is book-keeping, as we show that the existence
  of the sequence $\seq(x_{\alpha(r)})$ is incompatible with countable
  compactness. For each $r\in\nats$, let
  $A_r\defeq\set{x_{\alpha(s)}\st\text{$s\in\nats$, $s\gee r$}}$, and
  let $E_r\defeq\closure{A_r}$. Then $\seq(E_r)$~is a decreasing
  sequence of closed subsets of the countably compact space~$X$, so
  the intersection $E\defeq\Intersection_{r\in\nats}E_r$ is non-empty.
  Let $x$~be a point in~$E$. Since $F_i$~is closed we have $E\sub
  F_i$, and recalling that $F_i\sub U$ it follows that $x$~belongs
  to~$U_\beta$ for some $\beta<\omega_1$. Since $U_\beta$~is open, it
  is a neighbourhood of~$x$.

  There are now two cases to consider. Suppose first that
  $\beta\gee\alpha(r)$ for all $r\in\nats$. Since $\seq({\alpha(r)})$
  is increasing we must have $\beta>\alpha(r)$ for all $r\in\nats$,
  therefore $x_\beta\notin A_1$. Now $U_\beta\intersect
  Y=\set{x_\beta}$, and $A_1$~is a subset of~$Y$, so it follows that
  $U_\beta\intersect A_1=\emptyset$. But this contradicts the fact
  that $x$~belongs to~$E_1=\closure{A_1}$. Finally, suppose that
  $\beta<\alpha(r)$ for some natural number~$r$. Then we have
  $x_\beta\notin A_r$. Now $U_\beta\intersect Y=\set{x_\beta}$, and
  $A_r$~is a subset of~$Y$, so it follows that $U_\beta\intersect
  A_r=\emptyset$. Therefore $x\notin E_r=\closure{A_r}$, again a
  contradiction.
\end{proof}

The remaining result is my own combination of arguments from
\citet[3.12]{Juhasz80}, \citet[Prob.~III.68]{Arhangelskij84a} and
\citet[Lem.~4]{Arhangelskij71}.

\begin{prop}\index{countable tightness}\label{prop:regulartight-freesequence}
Let $X$~be a regular countably compact space of uncountable tightness. 
Then $X$~contains an uncountable subset~$Y$ which is discrete in the 
subspace topology.
\end{prop}

\begin{proof}
  The proof is split into two steps. In Step~2 we are going to
  construct an uncountable discrete subspace~$Y$ by transfinite
  induction. In order to facilitate this construction, we shall find
  in Step~1 a point~$x$ and a subset~$A$ such that $x$~belongs
  to~$\closure A\setminus A$, and $A$~meets every $G_\delta$~set that
  contains~$x$.
  \begin{enumerate}
  \item[Step~1] Since $X$~has uncountable tightness, there exists some
    non-closed subset of~$A$ that contains the closure of all its
    countable subsets. Pick any point~$x$ in $\closure A\setminus A$,
    and let $Q$~be a $G_\delta$~set that contains~$x$. We shall show
    that $Q$~meets~$A$ by the device of showing that $Q$~meets the
    closure of a countable subset of~$A$.

    Let $\set{U_i\st i\in\nats}$ be a countable family of open sets
    such that $Q=\Intersection_{i\in\nats}U_i$. Now, $X$~is a regular
    space, so for every $i\in\nats$ there exists an open set~$V_i$
    such that $x\in V_i\sub\closure{V_i}\sub U_i$. Let~$\Bee$ be the
    family of all finite intersections of $V_1$, $V_2$,~\dots, and let
    $P\defeq\Intersection_{i\in\nats}\closure{V_i}\sub Q$. Note that
    $\Bee$~is countable. Now, let $B=\Intersection_{i\in F}V_i$ be a
    typical element of~$\Bee$, where $F\sub\nats$ is finite. For each
    natural number~$i$, we have $x\in V_i$, therefore $x$~belongs
    to~$B$. As a finite intersection of open sets, $B$~is open, and
    since $x\in\closure A$ it follows that $B\intersect A$ is
    non-empty. For each $B\in\Bee$, let $x_B$~be a point
    in~$B\intersect A$, and let $C\defeq\set{x_B\st B\in\Bee}$. Then
    $C$~is a countable subset of~$A$, so therefore $\closure C\sub A$.
    We are going to show that $P\intersect\closure C\neq\emptyset$,
    and it will follow that $Q$~meets~$A$. Let $N$~be an open
    neighbourhood of~$P$; then $X\setminus N$ is a closed subset
    of~$X$, so it is countably compact. Now, recall that
    $P=\Intersection_{i\in\nats}\closure{V_i}$, therefore the family
    $\set{X\setminus\closure{V_i}\st i\in\nats}$ is an open cover for
    $X\setminus P$. The countably compact set $X\setminus N$ is
    contained in $X\setminus P$, so there exists a finite subfamily
    $\set{X\setminus\closure{V_1},\dots,X\setminus\closure{V_n}}$ such
    that $X\setminus N\sub\Union_{i=1}^n(X\setminus\closure{V_i})$.
    Therefore, letting $B\defeq\Intersection_{i=1}^nV_i$, we see that
    $B\sub\Intersection_{i=1}^n\closure{V_i}\sub N$. It follows that
    $x_B\in N$, so that $N\intersect
    C\contains\set{x_B}\neq\emptyset$. We can do this for every
    neighbourhood~$N$ of~$P$, hence $P\intersect\closure C$ is
    non-empty. This shows that $Q$~meets~$A$, which completes Step~1.

  \item[Step~2]\index{omega@$\omega_1$} We are going to define a
    family $Y=\set{y_\alpha\st\alpha<\omega_1}$ of distinct points
    of~$X$ such that~$Y$ is discrete as a subspace of~$X$. This will
    complete the proof of the proposition, because $Y$~is uncountable.

    Let $x$,~$A$ be as in Step~1. Pick any point~$y_{-\infty}$
    in~$A$, and let $W_{-\infty}\defeq X$. Let $\beta<\omega_1$, and
    suppose inductively that for all $\alpha<\beta$, we have found a
    point~$y_\alpha$ in~$A$ and an open set~$W_\alpha$ such that
    \begin{equation*}
      \text{(H$_\alpha$)}\;\left\{\;
        \begin{aligned}
          &x\in W_\alpha\text,\\
          &\closure{W_\alpha}\intersect
            \set{y_\gamma\st\gamma<\alpha}=\emptyset\text,\\
          &y_\alpha\in W_\gamma
          \qquad\text{for all $\gamma\lee\alpha$.}
        \end{aligned}
      \right.
    \end{equation*}
    Let $Y_\beta\defeq\set{y_\alpha\st\alpha<\beta}$. Then
    $Y_\beta$~is a countable subset of~$A$, so $\closure{Y_\beta}\sub
    A$. Therefore the point~$x\notin A$ from Step~1 does not belong
    to~$\closure{Y_\beta}$. So there exists an open neighbourhood~$N$
    of~$x$ such that $N\intersect Y_\beta=\emptyset$. Now $X$~is
    regular, so we can find an open set~$W_\beta$ such that $x\in
    W_\beta\sub\closure{W_\beta}\sub N$. Therefore
    $\closure{W_\beta}\intersect\closure{Y_\beta}=\emptyset$, that is
    $\closure{W_\beta}\intersect\set{y_\alpha\st
      \alpha<\beta}=\emptyset$. Now let
    $Q_\beta\defeq\Intersection\set{W_\alpha\st\alpha\lee\beta}$.
    Then $Q_\beta$~is a $G_\delta$~set that contains~$x$, so by Step~1
    above it has non-empty intersection with~$A$. Let~$y_\beta$ be
    any point in $Q_\beta\intersect A$. Then $y_\beta$~belongs to~$A$
    and we have $x\in W_\beta$ and $y_\beta\in W_\alpha$ for all
    $\alpha\lee\beta$. This shows that (H$_\beta$) is satisfied by
    our choice of $y_\beta$,~$W_\beta$. So by the principle of
    transfinite induction, we can define $y_\alpha$,~$W_\alpha$ in
    line with (H$_\alpha$) for all $\alpha<\omega_1$. Then
    $Y\defeq\set{y_\alpha\st\alpha<\omega_1}$ is an uncountable subset
    of~$X$. Now, by construction we have $y_\alpha\in W_\beta$ for
    all $\alpha\gee\beta$ but $y_\alpha\notin\closure{W_\beta}$
    whenever $\alpha<\beta$. It follows that, whenever
    $\alpha<\omega_1$, the set $Z\defeq
    W_\alpha\setminus\closure{W_{\alpha+1}}$ is open in~$X$ with the
    property that $Z\intersect Y=\set{x_\alpha}$. This shows that
    $Y$~is an uncountable discrete subspace of~$X$, as
    required.\qedhere
  \end{enumerate}
\end{proof}

Combining Propositions~\ref{prop:perfectcompact-countspread}--%
\ref{prop:regulartight-freesequence} we obtain
Theorem~\ref{thm:perfectcompact-counttight}.

One can also show that every compact perfectly normal space is
first countable. The proof will be omitted because it uses no new
ideas.

\subsection*{On Sequential Compactness}
\index{sequentially compact space|(} \index{countably compact
space|(}

When are compactness and sequential compactness equivalent? We saw
in Section~\ref{chap:sequential} that if the space is sequential, then
sequential compactness is equivalent to countable compactness.
Since compactness always implies countable compactness, it is then
natural to consider the question of when compactness is equivalent
to countable compactness. Trivially this is true in every Lindelöf
space. We also know that compactness and countable compactness are
equivalent in metric spaces; the following result is a
generalisation. The proof is based on that given in
\citet[Thm.~5.1.20]{Engelking89}.

\begin{lem}\index{paracompact space}
  \label{lem:paracompact+countcompact=compact}
  A regular countably compact paracompact space is compact.
\end{lem}

\begin{proof}
  This is true because all locally finite families in a countably
  compact space are in fact finite.  If a space $X$~is
  paracompact, then every open cover for~$X$ has a locally finite
  closed refinement, $\Eee$~say.  Suppose, to get a contradiction,
  that $\Eee$~is not finite.  Then $\Eee$~has a countable infinite
  subfamily $\set{E_i\st i\in\nats}$.  We shall show that $X$~is
  not countably compact. Since $\Eee$~is locally finite it is
  closure-preserving, so that $\Union\Dee$ is closed for every
  subfamily~$\Dee\sub\Eee$. Therefore if we let
  $F_1\defeq\Union_{i=1}^\infty E_i$, then $F_1$~is closed. Now
  let $F_2\defeq\Union_{i=2}^\infty E_i$. Then $F_2$~is also
  closed, and $F_2\sub F_1$. Continuing in this way we can define
  a decreasing sequence $\seq(F_i)$ of closed subsets of~$X$.
  Since $\Eee$~is locally finite, no point of~$X$ is contained in
  infinitely many elements of~$\Eee$, and it follows that no point
  of~$X$ is contained in all of the~$F_i$. Therefore
  $\Intersection_{i=1}^\infty F_i=\emptyset$. Thus $\set{F_i\st
    i\in\nats}$ is a countable family of closed subsets of~$X$ with
  the finite intersection property, but with empty intersection; this
  shows that $X$~is not countably compact.
\end{proof}

Of course this result may be strengthened. Indeed, the hypothesis
`paracompact' in the above can be replaced by `metacompact' as
shown by \citet[Thm.~2.4]{Arens50a}, or even by `pure' (see
\citet[Thm.~5]{Arhangelskij80}).

We may combine Lemma~\ref{lem:paracompact+countcompact=compact} and
Proposition~\ref{prop:sequential-compact}, to obtain a generalisation of
the result that compactness and sequential compactness are
equivalent in metric spaces.

\begin{cor}\index{paracompact space}
\label{cor:firstcount+paracompact+seqcompact=compact}
Let $X$~be a paracompact sequential space. Then $X$~is compact if and only 
if it is sequentially compact.\qed
\end{cor}

As stated earlier, there are reasons why one might speculate that 
compactness and countable compactness are equivalent for perfectly normal 
spaces. In fact, it is known that this cannot be decided within ZFC\@. 
\citet[p.~597]{Vaughan84} describes this as `one of the major results in 
set-theoretic topology.' 

%\pagebreak[3]
\begin{thm}\label{thm:perfectcountcompact-compact}\par\indent
\begin{enumerate}
\item 
Assume $\diamondsuit$, and let $X$~be Ostaszewski's space 
\citep[]{Ostaszewski76}. Then $X$~is perfectly normal, first countable and 
countably compact, but not compact.
\item 
Assume MA~and~$\neg$CH. Then every perfectly regular countably compact 
space is compact \textup(this result is due to 
\citet[Thm.~3]{Weiss78}\textup).
\end{enumerate}
\end{thm}

The whole question of when compactness and countable compactness are 
equivalent is considered in more detail in the survey by 
\citet{Vaughan84}. In particular the two parts of 
Theorem~\ref{thm:perfectcountcompact-compact} are there proved.

To get our final result on compactness in perfect spaces, we combine 
Theorem~\ref{thm:perfectcountcompact-compact} with 
Proposition~\ref{prop:sequential-compact}.

\begin{cor}
\label{cor:perfectfirstcount-compact}
The statement `Compactness and sequential compactness are equivalent
for perfectly regular sequential spaces' is independent of\/ ZFC\@.\qed
\end{cor}%
\index{countably compact space|)}%
\index{sequentially compact space|)}%
\index{compact space|)}%
\index{perfect space|)}

%%%%%%%%%%%%%%%%%%%%%%%%%%%%%%%%%%%%%%%%%%%%%%%%%%%%%%%%%%%%%%%%%%%%%
%%%%%%%%%%%%%%%%%%%%%%%%%%%%%%%%%%%%%%%%%%%%%%%%%%%%%%%%%%%%%%%%%%%%%
%%%%                                                             %%%%
%%%%                      BIBLIOGRAPHY                           %%%%
%%%%                                                             %%%%
%%%%%%%%%%%%%%%%%%%%%%%%%%%%%%%%%%%%%%%%%%%%%%%%%%%%%%%%%%%%%%%%%%%%%
%%%%%%%%%%%%%%%%%%%%%%%%%%%%%%%%%%%%%%%%%%%%%%%%%%%%%%%%%%%%%%%%%%%%%

%\pagebreak[4]

\end{document}